\documentclass[11pt]{amsart}
\usepackage{latexsym,amssymb,amsmath}
\usepackage{youngtab}
\def\BC{\mathbb C}\def\BF{\mathbb F}\def\BZ{\mathbb Z}
\def\BA{\mathbb A}
\def\BP{\mathbb P}
\def\pp#1{\mathbb P^{#1}}
\def\fa{\mathfrak a}
\def\fb{\mathfrak b}
\def\hh{\quad}
\def\fd{\mathfrak d}
\def\la{\lambda}
\def\ppp{{\mathbb P}}
\def\pp#1{{\mathbb P}^{#1}}
\def\tdim{\rm dim}
\def\hd{,...,}
\def\ww{\wedge}
\def\upperp{{}^\perp}
\def\TH{\Theta}

\def\inv{{}^{-1}}

\def\CC{\mathbb C}

\def\AA{{\mathbb A}}
\def\BB{{\mathbb B}}

\def\ZZ{\mathbb Z}

\def\11{\mathbf 1}
\def\PP{\mathbb P}

\def\fh{{\mathfrak h}}\def\fz{{\mathfrak z}}
\def\fs{{\mathfrak s}}
\def\fsl{{\mathfrak {sl}}}

\def\fso{{\mathfrak {so}}}
\def\fe{{\mathfrak e}}
\def\fu{{\mathfrak u}}
\def\ff{{\mathfrak f}}

\def\fg{{\mathfrak g}}
\def\fn{{\mathfrak n}}
\def\fp{{\mathfrak p}}

\def\fl{{\mathfrak l}}
\def\l{\lambda}
\def\a{\alpha}

\def\o{\omega}

\def\b{\beta}
\def\g{\gamma}
\def\s{\sigma}

\def\n{\nu}
\def\d{\delta}
\def\th{\theta}

\def\up#1{{}^{({#1})}}
\def\e{\varepsilon}
\def\ot{{\mathord{\,\otimes }\,}}
\def\op{{\mathord{\,\oplus }\,}}

\def\ra{{\mathord{\;\rightarrow\;}}}

\def\we{{\mathord{{\scriptstyle \wedge}}}}

\def\dim{{\rm dim}\;}
\def\La#1{\Lambda^{#1}}

\newcommand\rem{{\medskip\noindent {\em Remark}.}\hspace{2mm}}
\newcommand\exam   {{\medskip\noindent {\em Example}.}\hspace{2mm}}
\newcommand{\norm}[1]{\lVert#1\rVert}

\newtheorem{theo}{Theorem}[section]
\newtheorem{coro}[theo]{Corollary}
\newtheorem{lemm}[theo]{Lemma}
\newtheorem{prop}[theo]{Proposition}

\begin{document}

\title{Series of Lie groups}

\author{J.M. Landsberg and L. Manivel}

\begin{abstract}
For various   series  of complex semi-simple Lie algebras $\fg (t)$
equipped with irreducible representations $V(t)$, we decompose the tensor
powers of $V(t)$ into irreducible factors in a uniform manner, using
a tool we call {\it diagram induction}.
  In particular, we interpret  the
decompostion formulas of Deligne 
\cite{del} and Vogel \cite{vog} 
for decomposing
$\fg^{\ot k}$ respectively for the exceptional series and $k\leq 4$ and all
simple Lie algebras and $k\leq 3$, as well as
new formulas for the other rows of Freudenthal's magic chart. 
  By working with Lie algebras augmented by the symmetry group
of a marked Dynkin diagram, we are able to extend the list
\cite{brion}  
of modules for which the algebra of invariant regular functions
under a maximal nilpotent subalgebra is a polynomial algebra.
Diagram induction applied to the exterior algebra furnishes new examples
of distinct representations having the same Casimir eigenvalue.
\end{abstract}

\maketitle

\bigskip

\section{Introduction}

One way to define a collection of Lie algebras $\fg (t)$,
parametrized by $t$, each equipped with
a representation $ V(t)$, as forming
a \lq\lq series\rq\rq\ is, following Deligne, to require
 that the tensor powers of $V(t)$ should
decompose into irreducible $\fg (t)$-modules
in a   manner independent of $t$, with  
formulas for the dimensions of the irreducible components
of the form $P(t)/Q(t)$ with $P,Q$ polynomials decomposing into
products of  integral linear factors. 
In  this paper we study such decomposition formulas,  which provides a
companion to \cite{LM3} where we study the corresponding dimension
formulas. We connect the formulas to the geometry of the closed
orbits $X(t)\subset \BP V(t)$, and their unirulings by homogeneous
subvarieties. We relate the linear unirulings to work of Kostant
\cite{kos1}. By studying such series, we determine new modules
that, appropriately viewed, are {\it exceptional} in the sense of
Brion \cite{brion},  e.g., theorem 6.2.

The starting point of this paper was the work of 
Deligne et. al.  \cite{del, de2}, containing uniform
decomposition and dimension formulas for the tensor powers of the
adjoint representations of the exceptional simple Lie algebras 
up to $\fg^{\ot 5}$.
Deligne's  method for the decomposition formulas
was based on comparing Casimir eigenvalues and he offered a conjectural
explanation for the formulas via a categorical model based
on bordisms between finite sets.
 Vogel \cite{vog} obtained similar formulas for all simple Lie
super algebras based on his   {\it universal Lie algebra}.
We show that all  {\it primitive} factors in the decomposition formulas
of Deligne and Vogel
can be accounted for using diagram inductions.
(The non-primitive factors are those  either   inherited from lower
degrees or arising from a bilinear form, thus knowledge of
the primitive factors gives the full decomposition.) We also derive new decomposition
formulas for other series of Lie algebras.

In \S 2, we describe
a pictorial
procedure using Dynkin  diagrams for
determining the decomposition of $V^{\ot k}$. It was 
discovered by unifying several geometric
observations about the closed orbit $X=G/P\subset\BP V$.
We also give an interpretation of this diagram
induction in terms of vector bundles.

For example, in \cite{LM0} we 
  determined the Fano varieties  $\BF_k(X)$ of $\pp{k-1}$'s contained
in $X$. Since $\BF_k(X)$ is a subset
of the Grassmannian of $(k+1)$-planes through the
origin in $V$,  $G(k+1, V)\subset\ppp\La {k+1} V$, $\BF_k(X)$ determines a
distinguished component  of $\La {k+1}V$.
This component  has the property that  it
is a Casimir eigenspace and the corresponding Casimir
eigenvalue is maximal among Casimir eigenvalues of $\La {k+1}V$.
The components of $\La k\fg$ with maximal Casimir
eigenvalue were considered by Kostant in \cite{kos1}.  In \S 9 
we explain   how Kostant's results can be extended to 
general fundamental representations, 
 with special attention paid to minuscule representations.
 These Casimir eigenspaces provide new examples of distinct modules
 with the same Casimir eigenvalue which are different from the
examples in \cite{GKRS}.

In \S 3 and \S 4, we distinguish and interpret the {\it primitive}
components in the decomposition formulas of Deligne and Vogel.

The exceptional series of Lie algebras occurs as   a line in
 Freudenthal's magic square
(see e.g., \cite{freud}, or the variant we use in \cite{LM3}). 
The three other lines each come with   preferred representations.
Dimension 
formulas
for all representations supported on
the cone in the weight lattice
generated by the   weights 
of the preferred representations
similar to those of the exceptional series were obtained in \cite{LM3}. 
In \S 5, \S 6 and  \S 7
we obtain the companion decomposition formulas. 
A very nice property shared by  
many of these prefered representations 
is that they are {\it exceptional} in the sense 
of \cite{brion}, that is, their covariant algebras 
are polynomial algebras. 
We prove that in some cases   where this is not naively true, it becomes so 
if we take into account 
the symmetry group of the associated marked Dynkin diagram.

 In the course of revising the exposition of this paper for the
referee, we ran across the closely related preprint \cite{delgr}. 
  
\subsection{Notations and conventions}

For a given complex simple Lie algebra $\fg$ we fix a Cartan subalgebra
and set of positive roots.  
The highest root of $\fg$ (resp. 
the sum of the positive roots of $\fg$) will be denoted $\tilde\a$ (resp. 
$2\rho$).

Let 
$V=V_\la$ be be an irreducible representation 
 of highest weight $\l$ of $\fg$. To $V$ we associate
a marked Dynkin diagram
$ D(\fg ,\l)$   where
we identify the nodes of the diagram with the fundamental
weights $\o_1\hd \o_n$ and if $\l = a_1\o_1+\cdots + a_n\o_n$ we
mark the node corresponding to $\o_j$ with $a_j$.
We freely interchange
 the terminology \lq\lq marked Dynkin diagram\rq\rq\ and
\lq\lq irreducible representation\rq\rq .

Let $D(\ff)\subset D(\fg)$
be a subdiagram.
We define the {\it border set} of $D(\ff )$ in $D(\fg)$ to be 
the nodes of $D(\fg)\backslash D(\ff )$
adjacent to the nodes of $D(\ff)$. If $D(\fg, \l)$ is a marked
diagram where all the nonzero markings lie on $D(\ff)$,
we say $\l$ has {\it support} on $D(\ff )$ and   let
$W_{\l}$ denote the corresponding $\ff$-module.  

If $P$ denotes a  partition   of size $k$, we let 
$S_P(V)$ denote the corresponding
Schur power, which can be considered as a submodule of $V^{\ot k}$.

The Cartan product of 
two irreducible modules 
$V_{\mu}$ and $V_{\nu}$ is denoted  $V_{\mu}V_{\nu}:=V_{\mu+\nu} 
 \subset V_{\mu}\ot V_{\nu}$.

 For a given irreducible $\fg$-module $V_{\l}$, we let 
 $\th_{V_{\l}}=(\l,\l+2\rho)$ denote the Casimir
eigenvalue for $V$ with
the normalization $(\tilde\a,\tilde\a)=2$.

 We use the ordering of the roots as in \cite{bou}.

\section{Diagram induction} 

\begin{theo}
Let $\fg$ be a complex simple Lie algebra, with
Dynkin diagram $D(\fg)$ and let  $D(\ff )\subset D(\fg )$ be
a subdiagram, and let $\l$ be a weight of $\ff$ which we also consider as a
weight of $\fg$ as explained above. We let $V=V_{\l}$ (resp. $W=W_{\l}$)
denote the corresponding $\fg$ (resp. $\ff$-module).
Let $P$ be a partition of size $k$ and let $\o_1\hd \o_n$ denote
the fundamental weights of $\fg$.

1. If  
$\BC\subset S_P(W)$  (i.e., the Schur power $S_P(W)$ contains a trivial representation)  then the corresponding
Schur power $S_P(V)$ contains
an irreducible representation whose weight has 
support exactly
the border set $B$, i.e., the support is contained in $B$ and every
weight of $B$ appears with a nonzero multiplicity.

2. More generally if $W_{\eta}\subset S_P(W)$ is an irreducible submodule,
  write $\eta = k\l -\psi$, where $\psi$ is a sum of simple roots of
the root system of $\ff$. Consider $\psi$ as a sum of simple roots
of the root system of $\fg$, re-express $\psi$ as a sum of fundamental
weights of $\fg$, $\psi = a_1\o_1+\cdots + a_n\o_n$ and let $\tilde\eta=
k\l-(a_1\o_1+\cdots + a_n\o_n)$ denote the corresponding weight of $\fg$.
Then $\tilde \eta$ is a sum of weights from the border set    
of $\ff$ and the
weight $\eta$ considered as weights of $\fg$ and $V_{\tilde\eta}$ occurs
as a submodule of $S_P(V)$.

3. If $\l_1\hd \l_s$ are weights of $\ff$ and
$W_{\l_1}\ot \cdots \ot W_{\l_s}$ contains a trival representation, then
$V_{\l_1}\ot\cdots \ot V_{\l_s}$ contains
an irreducible representation whose weight has 
support exactly $B$ and the analogue of 2. holds for irreducible
submodules $W_{\eta}\subset W_{\l_1}\ot\cdots \ot W_{\l_s}$.
\end{theo}

  \medskip

\begin{exam} Let $(\fg, V_{\l})= (\fe_7, V_{\o_7})$,
let $(\ff, W_{\l})=(\fd_6, W_{\o_1})$. Then since
$S^2W$ contains the trivial representation, we have $V_{\o_7}\subset S^2V_{\o_1}$.

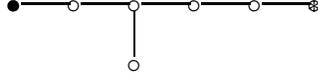
\begin{center}
\setlength{\unitlength}{4mm}
\begin{picture}(11,3)(0.8,-1.5)
\multiput(0,0)(2,0){6}{$\circ$}
\multiput(0.5,.3)(2,0){5}{\line(1,0){1.6}} \put(4,-2){$\circ$}
\put(4.25,-1.45){\line(0,1){1.5}}
\put(0,0){$\bullet$}
\put(10,0){$\ast$}
\end{picture} 
\end{center}

\end{exam}
\medskip

\begin{exam} Consider $V=V_{\o_k}=\La k\BC^n$
as a  $\fg=\fsl_n$-module.  
There is a natural quadratic form on
the $\fh=\fsl_{4p}$-module $\La  {2p}\CC^{4p}$. 
Thus for all $p<n$ the trivial representation
in $S^2\La{2p}\BC^{4p}$  induces a subspace of $S^2V$, namely
$V_{\o_{k-2p}+\o_{k+2p}}$ and these give us the full
decomposition
$$
S^2(\La  k \CC^n)
=S^2V_{\o_k}=V_{2\o_k}+ V_{\o_{k-2}+\o_{k+2}} +V_{\o_{k-4}+\o_{k+4}}
+\cdots
$$

Similarly, there is a   natural symplectic form on $\La  {2p+1}\CC^{4p+2}$ and
 the corresponding subdiagrams  recover the complete decomposition:

$$
\La  2(\La  k \CC^n)
=\La  2V_{\o_k}=  V_{\o_{k-1}+\o_{k+1}} +V_{\o_{k-3}+\o_{k+3}}
+\cdots
$$ 
\end{exam}

\begin{proof}
We have $\ff=\fs  + [\fs, \fs] $, where $\fs$ is   the sum of the 
root spaces   of $\fg$ with support in $D(\ff )$.  Every dominant 
weight $\s$ of $\fg$
with support in $D (\ff)$ defines a $\fg$-module $V_{\s}$ and an $\ff$-module
$W_{\s}$. 

The main observation is that 
if $\tau$ is a weight such that 
$\s-\tau$ has support on $D(\ff )$, its multiplicity must be  the same inside
$V_{\s}$ and $W_{\s}$. This is an easy  consequence of
Kostant's multiplicity formula:
$$\begin{array}{rcl}
\dim V_{\s}(\tau) & = & \sum_{w\in W(\fg )}\e (w) 
P(w(\s+\rho)-(\tau+\rho)), \\
\dim W_{\s}(\tau) & = & \sum_{w\in W(\ff)}\e (w)
 P'(w(\s+\rho')-(\tau+\rho')), 
\end{array}$$
where $P$ and $P'$ are the Kostant's partition functions in 
$\fg$ and $\ff$ respectively, and $\rho$, $\rho'$ the half
sums of the corresponding positive roots. First observe that 
since $\rho$ is also the sum of the fundamental weights,
if $w\notin W(\ff )$, then
  $w(\rho)-\rho$ will have some negative coefficient 
on a simple root not in $D(\ff )$. 
Then $\s-\tau+w(\rho)-\rho$, and a fortiori $w(\s+\rho)-(\tau+\rho)$, is not 
a sum of positive roots, and the partition function vanishes. 
Moreover, when $w\in W(\ff )$, we have $w(\rho')-\rho'=w(\rho)-\rho$, 
which proves   the two multiplicities  are the 
same. 

We give details for case 3 with two factors $W_\s\ot W_{\s'}$, the
other cases are similar:

Take two fundamental weights $\s$, $\s'$ with support in 
$D(\ff )$, and consider the decomposition  
$$W_{\s}\ot W_{\s'}=\bigoplus_{\tau\in \Xi}W_{\tau}.
$$ 
We can in principle obtain this decomposition by the 
following algorithm: we consider all the sums of a weight of
$W_{\s}$ with a weight of $W_{\s'}$ (with multiplicities); then
we choose a maximal element $\tau$ in this set. It must belong
to $\Xi$. Then we subtract the weights of $W_{\tau}$ (with their
multiplicities), and we continue until we are left with an empty
set of weights. 

Now consider the decomposition  of
$V_{\s}\ot V_{\s'}$. We   apply the   same procedure, only we
begin  with maximal weights which have support on $D(\ff )$. 
For these weights, the multiplicities are the same as
for the corresponding $\ff$-modules, thus we   obtain the same
set of dominant weights,
except that at the end we are left with
weights whose support is on $D(\fg )\backslash D(\ff)$
instead of the empty set. We conclude that
$$
V_{\s}\ot V_{\s'}\supset\bigoplus_{\tau\in \Xi}V_{\tilde\tau}.
$$ 
Here  
$\tilde\tau = \s +\s' -\psi_{\fg}$,
where $\psi_{\fg}$ is a non-negative linear 
combination of simple roots of $\fg$, is a weight of $\fg$
where the corresponding  weight $\tau$ of $\ff$ is
$\tau=\s +\s' -\psi$ where $\psi$ is the same non-negative linear 
combination of simple roots of $\fg$
only now considered as simple roots of $\ff$.

Comparing the expressions rewritten in terms of 
fundamental weights, we see that $\tau-\tau'$ is a linear 
combination of fundamental weights coming from $B$.
\end{proof}

\medskip

\subsection{Diagram induction via vector bundles}

We explain the case of inducing a representation from a trival
module $\BC\subset W_1\ot W_2$, the other cases being similar,
following the notation of above. (Here $W_1=W_{\s}$ etc... in the
notation above.)
Say $\BC$ induces a representation $U\subset V_1\ot V_2$. Let $\fp$
be the parabolic subalgebra of $\fg$ whose semi-simple Levi factor
is $\ff$. Pictorially,
$D(\ff )$ is
the subdiagram of $D(\fg)$ obtained by deleting the nodes
corresponding to $\fp$, with the  convention   that $\fp$ is generated by 
the root vectors 
corresponding to the Borel and  opposites  of the undeleted simple roots.

Consider the rational homogeneous variety $G/P\subset\ppp U$,
obtained by taking the projectivized orbit of a highest weight vector.
We  interpret diagram induction in terms of homogeneous
bundles on $G/P$. First of all, $U=\Gamma (L)$, i.e., $U$ is
the space of sections of
a homogeneous line bundle $L$ over $G/P$
and each $V_j$ is $\Gamma (E_j)$ for some homogeneous
vector bundle $E_j\ra G/P$.  
($L$ is the tautological (hyperplane) line bundle
on $G/P$.)

The homogenous vector bundles on
$G/P$ are in one to one correspondence with $P$ or $\frak p$-modules. Let
$W_j$ denote the irreducible $\fp$-module inducing $E_j$, i.e.,  $E_j=G\times_PW_j$.  
Write $\frak p = \ff \op \fz\op \fn$,
with $\ff$ semi-simple, $\fn$ nilpotent and $\fz$ the
center of the reductive part $\ff\op\fz$.

 $\fn$  acts trivially
on $W_j$ because $W_j$ is an irreducible $\fp$-module and
  $\fz$  acts by some character. 
We have a nonzero  multiplication map
$$
m: \Gamma(E_1)\ot\Gamma(E_2)\ra\Gamma (E_1\ot E_2).
$$
The occurence of $\BC$ in $W_1\ot W_2$ as an $\ff$-module
extends to a (nontrivial) $\fp$-submodule where $\ff$ and $\fn$ act  trivially
 and the new character for $\fz$ is the sum of the characters
  for $W_1$ and $W_2$,
thus we obtain a line subbundle $L\subset E_1\ot E_2$  and the desired inclusion
$U=\Gamma (L)\subseteq m(\Gamma(E_1)\ot \Gamma(E_2))$.

\begin{exam}
Consider $(\fg,V)=(\fsl (W),V_{\o_k}=\La k W)$. 
We have $U=V_{\o_{k-1}+\o_{k+1}}
\subset\La 2V_{\o_k} $ as mentioned above. Here $G/P=\BF_{k-1,k+1}$,
the variety of partial flags $\BC^{k-1}\subset\BC^{k+1}\subset W$,
$E$ is the   bundle whose
fiber over $(W_{k-1},W_{k+1})$ is
 $det(W/W_{k+1})\ot W_{k+1}/W_{k-1}$. (Due to our convention, we
 actually have $\Gamma (E)=V_{\o_k}^*$, not $V_{\o_k}$.)
Then $L= \La 2E$ has fiber $det(W/W_{k+1})\ot det(W/W_{k-1})$
and $\Gamma (\La 2 E)=V_{\o_{k-1}+\o_{k+1}}^*$.
\end{exam}

\medskip
We will apply diagram induction to study $\La kV$, $S^2V$ and
$S_{21}V={\rm Ker}(S^2V\ot V\ra S^3V)$. We first review some notions
of Tits:

\subsection{Tits' transforms and shadows}

For any  simple Lie group
$G$, with a fixed Borel subgroup, let $S,S'$ be two sets of positive roots of
$G$, and let $P_S$ be the parabolic subgroup generated by the Borel and the root 
subgroups generated by $-S$.   Consider the diagram 
$$\begin{array}{ccccc} & & G & & \\ &
{\scriptstyle \pi} \swarrow & &
\searrow {\scriptstyle \pi '} & \\
X= G/P_S & & & & X'= G/P_{S'}
\end{array}$$
Let $x'\in X'$ and consider $Y:= \pi (\pi'\inv (x'))\subset X$.
Tits calls $Y$ the {\it shadow} of $x'$ in $X$. Then $X$ is
covered by such shadows $Y$. Tits shows
\cite{tits} that $Y= H/R$ where
${\mathcal D}(H)={\mathcal D}(G)\backslash (S\backslash S')$, and 
$R\subset H$ is the parabolic subgroup corresponding to $S'\backslash S$.

\subsection{Submodules of $\La k V$}

We deduce the existence of  submodules of $\La kV_{\l}$ from marked subdiagrams of
the marked Dynkin diagram $ D(\fg , \l)$ isomorphic to $ D(\fa_{k-1} ,\o_1)$ 
via   the trivial representation $\La k\BC^k$.
 
In \cite{LM0} we showed that these subdiagrams describe linear unirulings of 
rational homogeneous varieties $X=G/P\subset\ppp V_{\l}$. 
If the subdiagram is of type $ D(\fa_{k-1} ,\o_1)$, we get a family of $\PP^{k-1}$'s 
on $X$ that are linearly embedded, hence  a linear uniruling of $X$. 
In the simply-laced case, all complete
families of linear unirulings are obtained that way. 

To recover a component of $\La k V$,  let $\BF_k(X)\subset G(k,V)\subset\PP(\La kV)$ 
be the Fano variety of $\pp{k-1}$'s in $X$  sitting inside
the Pl\"ucker embedded Grassmannian. Our uniruling defines a homogeneous component 
of $\BF_k(X)$, hence an irreducible submodule of $\La kV$ by taking the linear span.

In section \S 8  we determine the Casimir eigenvalues of these spaces. It
turns out that the linear span $\langle \BF_k(X)\rangle$ of the Fano variety 
is a Casimir eigenspace, with eigenvalue the largest possible.

For $k>1$, $\langle\BF_k(X)\rangle$ is usually not irreducible. When it is  reducible,
one obtains (pictorially!) distinct irreducible representations with
the same Casimir eigenvalue. This construction of representations
in the same Casimir eigenspace appears to be different from that
in \cite{GKRS}. Some of the Casimir eigenspaces $\langle\BF_k(X)\rangle$ were found
in \cite{WZ} via case by case checking. They were searching for such 
spaces because a homogeneous space $G/H$ with its  standard
homogeneous metric  is Einstein iff $T_{[e]}G/H$ is
a Casimir eigenspace of $H$.

\subsection{Submodules of $S^2V$}

We deduce the existence of  submodules of $S^2V_{\l}$ from marked 
subdiagrams  isomorphic to $D= D(\fd_{k } ,\o_1)$ or  $ D(\fb_{k } ,\o_1)$, 
thanks to the trivial representation given by the quadratic form.
 
\medskip In the case of such representations,
the highest weight $\tau$ of the induced
representation $V_{\tau}\subset S^2V$  can be
computed as follows (recall we already
have its support). Let $W$ be the Weyl group of $\fg$,
 and $W(D)\subset W$ the subgroup corresponding to the 
subdiagram $D$ so  $W(D)$ is generated by the 
simple reflections corresponding to the nodes of $D$. Let
$W_1(D)\subset W(D)$ be the stabilizer of $\lambda$, let $W^1(D)$ be the 
set of minimal length representatives of the cosets of $W_1(D)$ in
$W(D)$. Then $W^1(D)$ has a unique element $w_D^1$ of maximal length, 
and $\tau =\lambda +w_D^1(\lambda)$. 

We will use the notation $V_Q=V_{\tau}\subset S^2V$ to denote a representation
induced from a subdiagram of quadric type. By Tits'
transforms, the closed $G$-orbit $X_Q\subset\BP V_Q$ is a parameter space for
a uniruling of the closed orbit
 $X\subset \BP V$ by quadrics, i.e., linear sections $X\cap L$
that are quadric hypersurfaces in the projective space $L$. We use the
notation $Q=X\cap L$ to denote these quadrics. In the language of the
section above, these quadrics
are the shadows of points in $X_Q$ on $X$.


\begin{prop} Let $V=V_{\l}$ be a fundamental representation of $\fg$ such that
there is a subdiagram of quadric type $\fb_l$ or $\fd_l$ and let
$V_Q=V_{\tau}$ denote the induced submodule of $S^2V$. Then the 
Casimir eigenvalues are related by
$\th_{V_Q}=2(\th_V+ (\l,\l)-(\dim Q+2)(\a,\a))$, where $\a$ denotes the 
simple root dual to $\l$. 

In particular, if $\fg$ is simply-laced and $V=\fg$ is
the adjoint representation, then $\th_{\fg_Q}=2(\th_{\fg}-(\dim Q-1)(\a,\a))$.
\end{prop}

\begin{proof}We treat the case of $\fd_l$, the case of $\fb_l$ is similar.
Label the nodes of $D$ as $\a_1\hd\a_l$ and consider them as nodes
of $D(\fg)$ in what follows.
Let $\s = \a_1+\cdots +\a_{l-2}+\frac 12(\a_{l-1}+\a_l)$.
Note that with our normalizations,
 $(\s, 2\rho)=\tdim Q$,
$(\l,\s)=1$ and $(\s,\s)=1$. We have
$\tau = 2\l-2\s$ so $\th_{V_{\tau}}=(2\l-2\s, 2\l-2\s) +
(2\l-2\s, 2\rho)=2(\th_V + (\l,\l) - 4(\l,\s) + 2(\s,\s) -(\s,2\rho))$ and the result follows.
\end{proof}

Several such subdiagrams may exist, each of them will provide us 
with a component of $S^2V$.

\begin{exam} For every simple Lie algebra $\fg$ whose adjoint 
representation is fundamental, 
$S^2\fg$ contains only $\fg\up 2$, a trivial factor corresponding
to the Killing form, and factors of the form $\fg_Q$ (of which there
are at most three, or two up to a symmetry of the Dynkin diagram). 
\end{exam}

\begin{exam}In the case of the subexceptional (see \S 5) and Scorza series (see \S 6),
there is a unique $V_Q$ and  $S^2V=V\up 2 \op V_Q$. 
\end{exam}

\begin{exam}In the case of $(E_n, V_{\o_4})$ there are three distinct subdiagrams
of quadric type, but they furnish only a small part of $S^2V_{\o_4}$.
 \end{exam}
 
\medskip
 
  Note that in this case a point of $X_Q\subset \BP V_Q\subset \BP S^2V$
produces both a quadric hypersurface in $\BP V^*$ and a quadric section
of $X\subset \BP V$.

\smallskip
There is another characterization of 
maximal quadrics on $X=G/P\subset\PP V$, at least in the case
of minuscule and adjoint representations. Let  $\s_+(X)$ denote  
a component of the set of points 
of $\PP V\backslash X$ through which
 pass a family  of secants of $X$ of maximal dimension. 
If $p\in\s_+(X)$, its {\sl entry locus} $\Sigma_p=\{ x\in X, \exists y\in X-x, 
p\in\overline{xy}\}$, is a maximal quadric on $X$.  

\exam    Let $\fg=\fso_{2l}, V=V_{\o_k}$, with $1<k<l-1$. Here $X=G_Q(k,\BC^{2l})$
is the Grassmannian of $Q$-isotropic $k$-planes in $W=\BC^{2l}$, where $Q$ denotes 
the quadratic form  preserved by $\fg$.
The two families of quadrics given by Tits fibrations or diagram
induction  may be seen geometrically as follows: For the subdiagram
corresponding to $\fso_{2l-2k}$, choose $E\in G_Q(k-1,W)$, then
$$q_E=\{P\in G_Q(k,W), \; E\subset P\subset E^{\perp}\}\simeq Q^{2l-2k}.$$
The second family comes from the $\fa_3$ subdiagram. Pick
$E\in G_Q(k-2,W)$ and $F\in G_Q(k+2,W)$. Then
$$q_{E,F}=\{P\in G_Q(k,W), \; E\subset P\subset F\}\simeq Q^3.$$
We leave to the reader the pleasure of making the explicit correspondence
with the quadric hypersurfaces as above.  The correspondance with $\s_+(X)$ 
is straightforward: the line joining two distinct isotropic $k$-spaces 
$U, U'$ is contained in $X$ if and only if $U$ and $U'$ meet in codimension one, 
and $U+U'$ is isotropic. If this is not the case, points on  the secant 
line between $U$ and $U'$ are contained in $\s_+(X)$ if either $U, U'$ meet in 
codimension one but $U+U'$ is not isotropic
- in this case the entry locus is   $q_{U\cap U'}$,
unless $U, U'$ meet in codimension two and $U+U'$ is isotropic,
in which case the entry locus 
is   $q_{U\cap U',U+U'}$.

\subsection{Linear syzygies and subdiagrams}

Consider diagram induction when
$\ff=\fa_l$ with the
trivial representation in 
$W_{\tau_1}\ot W_{\tau_l}$. We    obtain
subrepresentations of $V_{\tau_1}\ot V_{\tau_l}$. We will call such
representations $(V_{\tau_1}V_{\tau_l})_{Aad}$.   
Changing notation, write 
$W_{\tau_1}=U_{\l}$, $W_{\tau_l}=W_{\mu}$, then
 $(UW)_{Aad}$  
has highest weight $\tau = \l+\mu - \s$ where $\s = \a_1+\cdots +\a_l$ where we
have labelled the roots corresponding to the subdiagram $D(\fa_l )$. 
We can thus compute its Casimir as above.

Let $S\subset S_{21}(V)$ denote the space of linear syzygies
among the generators of  $I(X)$, the ideal of $X$ (which are of degree two).
We have $S=S_{21}(V)\cap (I_2(X)\ot V)$. (We should really
consider $X\subset\ppp V^*$ here, but our abuse of notation is harmless.)
 
\begin{prop} Let
$V$ be a fundamental representation,
let $X\subset\ppp V$ be the closed orbit
and let  $U\subset I_2(X)$ be an irreducible component of the space of quadrics
containing $X$. Then $(VU)_{Aad}\subseteq S$.\end{prop}

Unfortunately we have no general proof of this fact, but it can be checked case by case. 

In the cases of the Severi and subexceptional
series below we   have equality.   

\medskip

\section{The Vogel decompositions}

Vogel \cite{vog} has proposed a {\it universal Lie algebra} $\fg ([\a,\b,\g ])$, 
which allows one to parametrize all complex
simple Lie superalgebras by a projective plane (over some 
extension of the rationals) quotiented out
by ${\mathfrak S}_3$. Evaluating at particular points, one recovers all
complex simple Lie algebras (and Lie superalgebras). He has given dimension and decomposition
formulae for the irreducible
modules in  $\fg^{\ot 2}, \fg^{\ot 3}$ that, independent of the existence of the
universal Lie algebra, give   decomposition and dimension formulae for actual
Lie algebras.

In order to connect his formulae to geometry, we   break the ${\mathfrak S}_3$ symmetry.
One reason for this is because inside $S^k\fg$ there is a preferred factor, the
Cartan power $\fg\up k$,
which has the geometric interpretation
of $I_k(X_{ad})\upperp$, the 
annhilator of the degree $k$ component of the ideal of the
closed orbit $X_{ad}\subset\BP (\fg)$.
 For example,  $\fg\up 2\subset S^2\fg$  could
be $Y_2, Y_2'$ or $Y_2''$
for Vogel (following his notations). We fix it to be $Y_2$. This has
the consequence of normalizing Vogel's parameter $\a$ to 
be $-(\tilde\a, \tilde\a)$, 
 as
according to Vogel, we have
$
2t=\th_{\fg}$
and
$
2(\th_{\fg}-\a)= 2\th_{\fg\up 2}= 2\th_{\fg}+2(\tilde\a,\tilde\a)
$.

In \S 2.3 we  discussed the factor $\fg_Q\subset S^2\fg$, where 
we take the largest subdiagram of quadric type here. We break the remaining
$\ZZ_2$ symmetry by requiring that this space be $Y_2'$. 
We obtain the following  
geometric interpretation of Vogel's parameter $\b$, which follows from 
\cite{vog}  and proposition 2.2:

\begin{prop} Notations as above: $\b= \tdim Q$ where $Q$ is the largest
  quadric contained in the adjoint variety $X\subset\ppp \fg$ obtained
  as a shadow as in \S 2.4.
\end{prop}

\exam    In the adjoint representations of $F_4$, $E_6$, $E_7$ and $E_8$, there is a unique 
quadric type subdiagram of the marked Dynkin diagram, respectively of types $B_3$, $D_4$, 
$D_5$ and  $D_7$, so that $\b=5,6, 8$ and $12$, respectively. 

\exam
If there is a second unextendable family of $G$-homogeneous quadrics on the
adjoint variety (as is the case for the orthogonal groups), then this supplies a geometric
interpretation of $\gamma$, namely $\gamma$ is the dimension of a quadric
in this second family. However, for adjoint representations,  this occurs only 
for the orthogonal groups, and in this case we always have $\gamma =4$. 

\medskip

Vogel describes three colinear collections of Lie algebras (in the sense
that some choice of inverse images of the
points associated with the Lie algebras are colinear in the projective plane).
The three Vogel lines are the exceptional, Osp, and Sl. To these
we add another line,
the {\it subexceptional series}, see \S 5, 
which lies on the line $2\a -\b+\gamma =0$.

With the
  above normalizations:
$$
\begin{array}{cccc}
 & \a &\b & \g \\
{\rm Exceptional} &  -2 & m+4 & 2m+4 \\
{\rm Osp:}\; SO(m), \;Sp(-m) & -2 & m-4 & 4 \\
{\rm Sl:}\; Sl(m) & -2 & 2 & m\\
{\rm subexceptional} & -2 & m & m+4
\end{array}$$

\medskip
In the exceptional series the values of $m$ are $-\frac{2}{3}, 0,
1, 2, 4, 8$ for $G_2,D_4, F_4, E_6, E_7, E_8$.  
The subexceptional line is
$A_1,   A_1\times A_1\times A_1,  C_3,  
  A_5,   D_6,   E_7$, with parameter
  $m=-\frac{2}{3},0,1,2,4,8$. Although
   $A_1\times A_1\times A_1$ is not simple,  
one can check that  the Vogel dimension and decomposition formulae still hold.
  The subexceptional line,
unlike the three other lines osp, sl and exceptional, is generic
to order three in the sense that none 
of the spaces that appear in Vogel's decomposition 
formulas are zero in $\fg^{\ot k}$ for $k\leq 3$
 except for the space Vogel labels $X_3''$ which is zero for
all simple Lie algebras. So, by comparing 
Casimir eigenvalues we can obtain geometric interpretations for all the Vogel spaces. These
interpretations (when such spaces exist)
 persist for other algebras not on the line.  
 
Here are Vogel's decompositions with our interpretations of the
spaces below. Recall our convention that $V_\mu V_\nu=V_{\mu+\nu}$.
$$
\begin{aligned}
\La 2 \fg &= X_1\op X_2\\
&=\fg\op\fg_2\\
S^2\fg &= Y_2\op Y_2'\op Y_2''\op X_0\\
&= 
\fg\up 2 \op \fg_Q \op \fg_{Q'} \op\BC_B\\
\La 3\fg&= 
X_3\op X_3'\op X_3'' \op X_2\op S^2\fg\\
&=\fg_3\op\fg_2\op S^2\fg\\
S^3\fg &= 2X_1\op X_2\op B\op B'\op B''\op Y_3\op Y_3'\op Y_3''\\
&=
 2\fg\op\fg_2\op B\op\fg\fg_Q\op\fg\fg_{Q'}
\op\fg\up 3\op\fg_{\BA\pp 2}\op Y''\\
S_{21}\fg &= 2X_2\op 2X_2\op Y_2\op Y_2'\op Y_2''\op B\op B'\op B'' \op C\op C'\op C''
\\
&= 2\fg\op2\fg_2\op\fg\up 2\op\fg_Q\op\fg_{Q'}\op B\op \fg\fg_Q\op\fg\fg_{Q'}
\op\fg\fg_2\op
(\fg\fg_Q)_{Aad}\op (\fg\fg_{Q'})_{Aad}
\end{aligned}
$$

We have written $\fg_3=X_3\op X_3'$
as it is a Casimir eigenspace.
We have no interpretation for $X_3''$ as it does not exist for
actual Lie algebras, nor $B$ because it does not exist for the exceptional series
and it is $-\fg_2$ for the subexceptional series.

The other decompostions can be deduced from these,
e.g., $X_1\ot X_2=\fg\ot\La 2\fg - \fg\ot\fg=(S_{21}\fg\op\La 3\fg) -(S^2\fg
\op\La 2\fg)$. 

The only space not yet explained is $\fg_{\BA\pp 2}$. It comes from
diagram induction applied to a subdiagram in the {\it Severi series},
the distinguished representations in the second row of Freudenthal's
magic chart (see \S 7)
as there is an invariant cubic on the representation $W$.
The subdiagram for $\fg_{\BA\pp 2}$  in the exceptional line is obtained 
by deleting the nodes for $\fg$ and $\fg_2$.

\subsection{Comparison with Freudenthal's magic square}

Normalizing $\a=-2$, Vogel's formula for $\tdim\fg$ is
$$
\tdim \fg =\frac{(\b+\g-1)(2\b+\g-4)(2\g+\b-4)}{\b\g}
$$

The triality model enables one to deduce the
following two parameter formula for the dimensions of
the Lie algebras occuring in Freudenthal's magic square (see \cite{LM3}).

\begin{prop}   Let $\fg(a,b)$ denote the 
Lie algebra that Freudenthal associates to the pair $(\AA,\BB)$ of
real division algebras of dimensions $a$ and $b$,
and let $p=a+4,q=b+4$ then 
$$\dim\fg(a,b)=3 \frac{(ab+4a +4b-4)(ab+2a+2b)}{(a+4)(b+4)}.$$
$$\dim\fg(p,q)=
\frac{3 (p q-20) (p q-2 p-2 q)}{pq}
$$
\end{prop}

The $a,b$ parametrization is natural from the point
of view of the composition algebras, the $p,q$ parametrization is
more natural from the point of view of Tit's fibrations.
That the $p,q$ parametrization might be simpler to work with   was brought to
our attention by B. Westbury.

\section{The  exceptional series}

Using either Freudenthal's perspective of incidence geometry \cite{freud} or
the triality model \cite{LM3}, one has four distinguished
representations in the exceptional series, denoted $X_1, X_2, X_3, Y_2^*$
in \cite{del}. We refer the
reader to \cite{com} for the notations and decomposition 
formulae.

The spaces $Y_3^*, G^*, H^*, I^*, Y_4^*$ all contain virtual
representations, i.e., negatives of actual representations,
  for the
larger algebras in the series
so it is not possible to assign direct geometric
interpretations.  

The primitive representations are as follows:
$X_k=\fg_k, Y_2^*= \fg_Q, C^*=(\fg I_2)_{Aad}=S_1,
 F^*=(C^*\fg)_{Aad}\subseteq S_2$. Here $S_1,S_2$ denote
 the first and second linear
 syzygies among the quadrics in the ideal for
 the closed orbit $X_{ad}\subset \BP \fg$,  where in 
 general, for an algebraic variety
 $X\subset \BP V$ defined by quadratic polynomials $I_2(X)\subset
 S^2V^*$,  we let $S_1:= (V^*\ot I_2(X))\cap S_{21}V^*$    and 
 $S_2:= (V^*\ot S_1)\cap S_{211}V^*$ the second 
 linear syzygies.

The others can be deduced from the primitive ones through Cartan products:
$Y_k=\fg\up k, A=\fg Y_2^*, C=\fg\fg_2, D=Y_2^*\fg\up 2, D^*=Y_3^*\fg,
E=\fg C^*,  F=\fg_2Y_2^*, G=\fg_2\fg\up 2, H=\fg_2\up 2, I=\fg\fg_3,
J={Y_2^*}\up 2$.

\medskip   
 
\section{Subminuscule representations}

Recall that a $\fg$-module $V$ is called of {\it type-$\theta$} if there is
a $\BZ/m\BZ$-grading  (allowing the possibility of $\BZ$-gradings
as well)
of a simple Lie algebra $\fl$ such that $\fg$ is
the semi-simple part of $\fl_0$ and $V=\fl_k$ for some $k$. It is
of {\it type I-$\theta$} if $k=1$ and the grading is a $\BZ$-grading. We
define a further subclass, the {\it sub-minuscule} representations
where the grading of $\fl$ is minuscule (i.e., three step). Geometrically,
the subminuscule representations are the representations of semi-simple Lie
algebras occuring as the isotropy representation on the tangent
space of an irreducible compact Hermitian symmetric space, the
type I-$\theta$ representations occur as the submodules $T_1\subset T_{[e]}G/P$
where $P$ is a maximal parabolic and $T_1$ is the (unique)
irreducible $P$-submodule of $T_{[e]}G/P$ see \cite{LM0}.

In \cite{LM0}
we showed that for subminuscule representations,
 the only $G$ orbits in $\ppp V$ are the
smooth points of the successive secant varieties of the closed
orbit $X=G/P\subset\ppp V$, and moreover that the union
of the secant $\pp{k-1}$'s, denoted $\s_k(X)$, is such that its
ideal is generated in degree $k+1$ with
$I_{k+1}\s_k(X)=I_2(X)\up{k-1}:=(I_2(X)\ot S^{k-1}V^*)\cap S^{k+1}V^*$,
where $I_2(X)\up{k-1}$ is called the $(k-1)$-st
{\it prolongation} of $I_2(X)$. 
Another way to phrase this prolongation property is
that the spaces of generators are the sucessive Jacobian
ideals of the highest degree space of generators.
 In practice these spaces are quite
easy to compute and
comparing with \cite{brion},
we observe that the symmetric algebra 
is free and the prolongations
of $I_2(X)$   furnish all the primitive factors for the
symmetric algebra,  so we obtain:

\begin{theo} Let $V$ be a sub-minuscule representation of
a semi-simple Lie algebra $\fg$.
With the notations above, and our convention
$V_{\mu}V_{\sigma}=V_{\mu +\sigma}$, we have a uniform formula
for the decomposition of the symmetric algebra into irreducible
$\fg$-modules:
$$
\oplus_{k=1}^{\infty} t^kS^kV= \Pi_{j=2}^{\infty} (1-t^jI_2(X)\up{j-1})^{-1}
$$
The product on the right hand side is finite.
\end{theo}

 Here the orbit closures exactly provide the primitive
factors for the symmetric algebra. In general, the orbit closures
will provide some, but not all primitive factors, see the examples of the subexceptional
and sub-Severi series below.

\rem A version of this result appears to have been known to
Kostant as the \lq\lq cascade of orthogonal vectors\rq\rq .

\subsubsection*{Example: The Scorza series}
Zak established an upper bound
  on the codimension of a smooth variety
$X^n\subset\pp{n+a}$
  of a given secant defect. (The secant defect is the
  difference between the 
  expected dimension  of the secant variety of $X$
  (min$\{n+a,2n+1\}$) and its 
  actual dimension.)  He then went on
to classify the varieties achieving this bound, which he calls the
{\it Scorza varieties}.
They are all closed orbits $G/P\subset\ppp V$
and give rise to the following two parameter $(a,n)$ series: $(SL_n, V_{2\o_1})$, 
$(SL_n\times SL_n, V_{\o_1}\ot W_{\eta_1})$,
$(SL_{2n}, V_{\o_2}), (E_6, V_{\o_1})$, where $a$ is respectively $1,2,4,8$
and $a=8$ only for the $n=3$ case.
 This series is the second row of the  generalized Freudenthal magic square,
 see \cite{LM1}.  
We could add to this the finite group ${\mathfrak S}_n$, case $a=0$, corresponding
to the variety of $n$ points in $\pp{n-1}$. 
  In this case the symmetric algebra is
generated by the \lq\lq determinant\rq\rq\
(see \cite{LM2}), which has  degree $n$, and the spaces of $k\times k$
minors. Here  $I_2(X)\up{j-1}$ respectively  has highest weights
$2\o_{n-k}$, $\o_{n-k}+\eta_{n-k}$, $\o_{2n-2k}$.
We remark that $\dim V(a,n)= n+a \frac{n(n-1)}{2}$.

\section{The subexceptional series}

This is the series coming from the third line of Freudenthal's 
square:
$$A_1, \quad A_1\times A_1\times A_1,\quad C_3,  
\quad A_5, \quad D_6, \quad E_7.$$

Let $m=-\frac{2}{3},0,1,2,4,8$ respectively.
Freudenthal's perspective \cite{freud} or the triality model
\cite{LM3} uncovers
three preferred irreducible
representations, respectively denoted $V, V_Q=\fg, V_2$
in the table below and of dimensions
$6m+8, \frac{3(2m+3)(3m+4)}{(m+4)}, 9(m+1)(2m+3)$.


Let $\Gamma_0$ be the
automorphism group of the Dynkin diagram, and $\Gamma\subset\Gamma_0$ 
be the subgroup preserving the marked Dynkin diagram 
($\Gamma={\mathfrak S}_3$ for $\fg=A_1\times A_1\times A_1$,
$\Gamma={\mathfrak S}_2$ for $\fg=A_5$ and is trivial otherwise).
   With the help of the programm LiE
\cite{LiE}, we obtained the following 
decomposition formulae into   $\fg\times\Gamma_0$-Casimir eigenspaces.
Letting $V_0=\CC  $, we have, up to at least degree six:
$$\begin{array}{rcl}
\La {2p}V & = & V_{2p}\op  V_{2p-2}\op\cdots\op V_0, \\
\La {2p+1}V & = & V_{2p+1}\op  V_{2p-1}\op\cdots\op V_1.
\end{array}$$
Note that $V_2,V_3$ are irreducible.

This decomposition coincides with the decomposition into
primitives for the symplectic form. In general the decomposition
of a symplectic $\fg$ module $W$  into primitives
is not Casimir-irreducible. Consider the primitives in the
$A_9$-module $\Lambda^2(\Lambda^5\BC^{10})=\BC \op V_{\o_2+\o_8}
\op V_{\o_4+\o_6}$. The last two factors consitute the primitive subspace
but they have different Casimir eigenvalues.

\smallskip
  In the last 
four cases of the series, $V$ is {\it exceptional} in the sense of \cite{brion}, 
that is, the algebra $\CC[V]^{\fu}$ of  
invariant regular functions on $V^*$ under a maximal
 nilpotent subalgebra $\fu$ of $\fg$, i.e., the covariant algebra,  is 
a polynomial algebra. Such an invariant is
   a highest weight vector of some 
symmetric power of $V$, which   allows 
one to decompose $S^kV$   into 
irrreducible factors
for all $k$. The results of \cite{brion} imply, 
again with our convention $V_{\l}V_{\mu}=V_{\l+\mu}$, that:
$$\bigoplus_{k\ge 0}t^k S^kV = 
(1-tV)^{-1}(1-t^2\fg)^{-1}(1-t^3V)^{-1}(1-t^4)^{-1}(1-t^4V_2)^{-1}.$$

As with the subminuscule case, 
the spans of generators of ideals
of each orbit closure in $\BP V$ give primitive factors in 
$S^\bullet V$. In contrast, there is one additional primitive factor,
$V_2\subset S^4V$, which is also  the primitive part of $\La 2 V$. 
The presence of the primitive $V_2$ factor may be understood
as follows: the symplectic form
$\o$ on $V$ enables an equivariant identification $V\simeq V^*$. 
Polarizing  the invariant quartic form  
gives a map $q : S^3V\ra V^*\simeq V$. Finally, we obtain a natural map 
$s : S^4V\ra V_2$ by letting $s(v^4)=v\we q(v^3)$ mod $\o$.
 This map is nonzero and exhibits $V_2$ as an 
 irreducible component of $S^4V$. \medskip

\medskip
 The remaining decompositions for $V^{\ot k}$ in degrees three and four are:
 $$\begin{array}{rcl}
 S_{21} V &= & V\op C\op {V\fg}\op {VV_2}\\
 S_{31} V & =& {V_2}\op 2{V\up 2}\op 2\fg\op {VC}\op {\fg V\up
2}\op {\fg V_2}\op {\fg_2}\op {V_2V\up 2}, \\
S_{22}V & =& \CC\op 2{V_2}\op {\fg V\up 2}\op {VC}\op Q\op
{\fg\up 2}\op {V_2\up 2}, \\
S_{211}V & =& {V_2}\op {V\up 2}\op\fg\op {VC}\op {\fg
V_2}\op {\fg_2}\op L\op {VV_3},
\end{array}
$$

Note   that the only primitives
up to degree three are $C, \fg$ and the $V_k$'s and the only new
primitives in degree four are $Q$ and $L$.

The   Casimir eigenvalues for the modules involved in 
  these formulas are all of the form $\frac{am+b}{8m+8}$ with $a,b\in\BZ$,
  and are 
linear functions of $(\la,\la )  =\frac{3}{8m+8}$. Here are the Casimir eigenvalues:
$$\begin{array}{ccc}
\th_{V\up k}=\frac{k(6m+9)+3(k^2-k)}{8m+8},
&  \th_{V_k}=\frac{6km+ (10k-k^2)}{8m+8},
& \th_{\fg\up k}=\frac{2km+ ( k^2+k)}{2(m+1)},  \\
  \th_C=\frac{12m+9}{8m+8}, & 
  \th_Q=\frac{3m}{2m+2}, & 
 \th_L=\frac{2m+1}{m+1}
\end{array}$$

The dimensions of these modules are rational functions of $m$
with simple denominators, see \cite{LM3} for the dimension formulas
with the exception of
$$\dim C = \frac{32(m+1)(2m+3)(3m+4)}{(m+4)(m+6)}, $$
$$\dim Q \; = \; \frac{(8-m)(m+1)(2m+3)(3m+2)(3m+4)}{(m+4)^2(m+6)}, $$
$$\dim L \; = \; \frac{9(8-m)(m+1)(2m+3)(3m+2)(3m+4)}{(m+4)(m+6)(m+8)}. $$

There is a geometric interpretation for the primitives $C$ and $L$ in terms of syzygies.
 We lack a geometric interpretation for $Q$ as it
is empty for $\fe_7$ and it does not appear in the minimal free resolutions. 

\begin{prop} 
Let  $S_k$ denote the space of linear syzygies of order $k$  
in the minimal resolution of a sub-exceptional variety, beginning 
with $S_0=I_2(X)$. Then
$$S_0=\fg, \qquad S_1=C, \qquad S_2=L.$$
\end{prop}

The   decompositions of $\fg^{\ot 2}$ and $\fg^{\ot 3}$ are as with Vogel's formulas. 
Except in the case of $\fe_7$, where
it is irreducible,   $\fg_3 $ decomposes
into two irreducible representations
that are called  $X_3$ and $X'_3$, by Vogel.  Their
dimensions  have algebraic expressions which are not rational in $m$.
(In Vogel's formulae, the expressions are not rational in
$\a,\b,\g$ either.)  From Deligne's perspective,
$\fg_3$ should not be considered a preferred representation as
its dimension formula contains a quadratic factor in its numerator:
$$\dim \fg_3= \frac{ (2m+3)(3m+4)(9m+16)(m+1)(18m^2+43m+4)}{(m+4)^3}
$$

\medskip We also have:

$$\begin{array}{rcl}
\fg\ot V & = & V\op C\op {V\fg},  \\
{\fg\up 2}\ot V & =& {V\fg\up 2}
\op {V\fg}
\op {\fg C}, \\ {V\up 3}\ot\fg & =& {\fg V\up 3}\op {V\up 3}\op
{VV_2}\op {CV^2},
\\
\La 2{V_2} & = & {V\up 2}\op \fg\op {VC}\op {\fg V_2}\op
{\fg_2}\op L\op {VV_3},
\\
\La 2{V\up 2} & = & {V\up 2}\op \fg\op {\fg V_2}\op
{V_2V\up 2}, \\
S^2{V\up 2} & = & \CC\op {V_2}\op {\fg V\up 2}\op {V\up
4}\op \fg\up 2\op {V_2\up 2}, \\ {V_2}\ot {V\up 2} & = & {V_2}\op {V\up
2}\op
\fg\op {VC}\op {\fg V\up 2}\op {\fg V_2}\op {\fg_2}\op {VV_3}\op {V_2V\up
2},
\\ {V_2}\ot\fg & = & {V_2}\op {V\up 2}\op \fg\op {VC}\op {\fg V_2}\op
{\fg_2}\op L, \\ {V\up 2}\ot\fg & = & {V_2}\op {V\up 2}\op {VC}\op {\fg
V\up 2}, \\ 
{V_2}\ot V & = & C\op {V\fg}\op {V_3}\op {VV_2}\op {V}, \\
{V\up 2}\ot V & = & {V\fg}\op {VV_2}\op V\up 3\op {VV_3}, \\
C\ot V & = & {V_2}\op\fg\op {VC}\op {\fg_2}\op L\op Q, \\
{V\fg}\ot V & = & {V_2}\op {V\up 2}\op \fg\op {VC}\op {\fg V\up 2}\op
{\fg V_2}\op {\fg_2}\op \fg\up 2,
\\ {V_3}\ot V & = & {V_2}\op {VC}\op {\fg V_2}\op {\fg_2}\op L\op
{VV_3}\op V_4, \\ {VV_2}\ot V & = & {V_2}\op {V\up 2}\op {VC}\op {\fg
V\up 2}\op {\fg V_2}\op {VV_3}\op {V_2V\up 2}\op {V_2\up 2}, \\ V\up 3\ot V
& = & {V\up 2}\op {\fg V\up 2}\op {V_2V\up 2}\op {V\up 4}, \\
S^2 V_2& = &V_2\up 2 \op V_4 \op \fg V\up 2 \op \fg\up 2 \op
2V_2
\op \CC \op CV
\end{array}$$

\medskip
The highest weights of the modules involved in the above formulas
are as follows:

$$\begin{array}{lrrllll} 
 & A_1 & A_1^{\op 3} & C_3 & A_5 & D_6 & E_7 \\ 
  & & & & & & \\
V & \left[3\right] & \left[1,1,1\right] & \left[0,0,1\right] &
\left[0,0,1,0,0\right]
& \left[0,0,0,0,0,1\right] & \left[0,0,0,0,0,0,1\right] \\
V_2 & \left[4\right] & \left[2,2,0\right] & \left[0,2,0\right] &
\left[0,1,0,1,0\right]
& \left[0,0,0,1,0,0\right] & \left[0,0,0,0,0,1,0\right] \\
{V_3} & &\left[3,1,1\right] & \left[1,2,0\right] & \left[1,0,0,2,0\right]
& \left[0,0,1,0,1,0\right] & \left[0,0,0,0,1,0,0\right] \\
{V_4} &-\left[ 4\right] & [4,0,0] & \left[0,3,0\right] & \left[0,3,0,0,0\right]
& \left[0,1,0,0,2,0\right] & \left[0,0,0,1,0,0,0\right] \\
  & & [2,2,2] &  \left[3,0,1\right]  &   \left[1,1,0,1,1\right]
& \left[0,0,2,0,0,0\right] &  \\
{V_5} &-\left[ 3\right] & & \left[ 2,1,1\right]  &
\left[2,0,1,0,2\right]   & \left[1,0,0,0,3,0\right] &
\left[0,1,1,0,0,0,0\right] \\
  &  & & \left[ 5,0,0\right]  & \left[ 1,2,0,1,0\right]
& \left[0,1,1,0,1,0\right]  &   \\
{V_6} &-\left[ 0\right] & -[4,0,0] & \left[ 4,1,0\right]   &  \left[
2,1,1,0,1\right] & \left[0,0,0,0,4,0\right] & \left[1,2,0,0,0,0,0\right]\\
  & & -[2,2,2]& \left[ 2,0,2\right]  & \left[ 0,2,0,2,0\right]  &
  \left[1,0,1,0,2,0\right]  & \left[0,0,2,0,0,0,0\right]\\
  & & &  & \left[ 3,0,0,0,3\right]
&  \left[0,2,0,1,0,0\right] &  \\
\fg & \left[2\right] & \left[2,0,0\right] & \left[2,0,0\right] &
\left[1,0,0,0,1\right]
& \left[0,1,0,0,0,0\right] & \left[1,0,0,0,0,0,0\right] \\
{\fg_2} & &\left[2,2,0\right] & \left[2,1,0\right] &
\left[0,1,0,0,2\right]  & \left[1,0,1,0,0,0\right] &
\left[0,0,1,0,0,0,0\right] \\
{\fg_3} & & [2,2,2] & \left[3,0,1\right] &
\left[3,0,1,0,0\right]  & \left[2,0,0,1,0,0\right] &
\left[0,0,0,1,0,0,0\right] \\
  & & -[4,0,0] & \left[ 0,3,0\right] &
\left[1,1,0,1,1\right]  & \left[0,0,2,0,0,0\right] & \\
C & \left[1\right] & \left[1,1,1\right]\ot\rho & \left[1,1,0\right] &
\left[1,1,0,0,0\right]
& \left[1,0,0,0,1,0\right] & \left[0,1,0,0,0,0,0\right]\\
Q & & \left[0,0,0\right]\ot\rho & \left[0,1,0\right] &
 \left[1,0,0,0,1\right]
& \left[2,0,0,0,0,0\right] & \\
 L & &
\left[2,0,0\right] &
\left[1,0,1\right] & \left[0,1,0,1,0\right]  &
\left[0,0,0,0,2,0\right] & 
 \end{array}$$

\bigskip In the column corresponding to $A_1\times A_1\times A_1$, 
$\rho$ denotes the two-dimensional irreducible representation of 
$\Gamma={\mathfrak S}_3$. 
\medskip

The first two cases of the series deserve special care since they are slightly 
degenerate and we discuss them in
 the following two subsections. 

\subsection
{Binary cubics}   
 In the $A_1$ case $V_2= \fg\up 2$, and there is no factor $1-t^4V_2$ in the denominator. 
Moreover, $V$ is not exceptional 
since there exists a relation in degree $6$ between the fundamental covariants 
(see e.g. \cite{dixm}). We have  
$$\bigoplus_{k\ge 0}t^k S^kV = \frac{1-t^6V^{(2)}}
{(1-tV)(1-t^2\fg)(1-t^3V)(1-t^4)}.$$

\subsection
{ 
$2\times 2\times 2$ hypermatrices} 

Write $A_1\times A_1\times A_1=\fsl (A) \times \fsl (B)\times \fsl (C)$,
with $A,B,C\simeq \BC^2$.

Introduce the symmetrization operator $\phi$  on formal power series with 
coefficients in $A_1\times A_1\times A_1$-modules, which associates to 
$S^aA\ot S^bB\ot S^cC$ its complete symmetrization, e.g., $\phi(S^2A\ot B)=S^2A\ot B\op
S^2A\ot C\op S^2B\ot A\op S^2B\ot C\op S^2C\ot A\op S^2C\ot B$, and $\phi(S^3A)=
S^3A\op S^3B\op S^3C$. 

\begin{theo}The covariant 
algebra $\CC[A\ot B\ot C]^{\fn\times {\mathfrak S}_3}$ is a polynomial algebra.
More precisely,
$$\bigoplus_{k\ge 0}t^k S^kV = \phi\frac{1}{(1-tV)(1-t^2\fg)
(1-t^3V)(1-t^4)(1-t^4V_2)},$$
where $V=A\ot B\ot C$,  $\fg =S^2A\op S^2B\op S^2C$ 
and $V_2=\La 2 (A\ot B\ot C)/\CC = S^2A\ot S^2B\op 
S^2B\ot S^2C\op$ \linebreak $ S^2C\ot S^2A$.
\end{theo}

 Here use the convention $\fg^{(k)}= S^{2k}A\op S^{2k}B\op S^{2k}C$,
  and   similarly for $V_2^{(k)}$.

Thus although $A\ot B\ot C$ is not exceptional in the sense of \cite{brion}, 
it does become exceptional when we take into 
account the ${\mathfrak S}_3$-symmetry.

Note that the  generators of the symmetric algebra have
  the same degrees as in the other cases of the subexceptional series. \medskip

The theorem is a consequence of the following lemma:

\begin{lemm} 
Let $\mu(n;a,b,c)$ denote the multiplicity of $S_{n-a,a}A\ot S_{n-b,b}B\ot S_{n-c,c}C$ inside 
$S^n(A\ot B\ot C)$. Suppose that $c\ge a,b$ and $2c\le n$. Then 
$$\mu(n;a,b,c) = \left\{\begin{array}{lcl}
0 & {\rm if} & c>a+b, \\ E(\frac{a+b-c}{2})+1 & {\rm if} & c\le a+b
 \;{\rm and}\; n\ge a+b+c, \\
E(\frac{a+b-c}{2})-E^+(\frac{a+b+c-n}{2})+1 & {\rm if} & c\le a+b \;{\rm and}\; n\le a+b+c.
\end{array} \right .$$
\end{lemm}
 
Here $E(x)$ denotes the largest integer smaller than or equal to $x$, 
and $E^+(x)$ the smallest integer greater than or equal to $x$. 

 Recall that irreducible 
representations of ${\mathfrak S}_n$ are naturally indexed by partitions
 of $n$. We let   
$[\l]$ denote the representation associated to a partition $\l$, following the
notation of \cite{man2}.  

By Schur duality, $\mu(n;a,b,c)$ can be 
interpreted in terms of representations of symmetric groups, as the dimension of the 
space of ${\mathfrak S}_n$-invariants in the triple tensor product $[n-a,a]\ot [n-b,b]\ot [n-c,c]$, 
or the multiplicity of $[n-a,a]$ inside $[n-b,b]\ot [n-c,c]$.  The behavior of the multiplicity 
of $[n+\l]$ inside $[n+\mu]\ot [n+\n]$ as a function of $n$ was investigated 
in \cite{brion3,
man2}, where it was proved to be non-decreasing, and constant 
for $n$ sufficiently large. 
\medskip

\proof We use Cauchy formula \cite{litt}  for the symmetric powers of a tensor product:
$$\bigoplus_{k\ge 0}t^k S^k(A\ot B\ot C) = \bigoplus_{a\ge b\ge 0}
t^{a+b}S_{a,b}A\ot S_{a,b}(B\ot C).$$ 
Since $A$ is two dimensional,   $S_{a,b}A=S^{a-b}A$ as $\fsl_2$-modules. 
Moreover, we can write 
$S_{a,b}(B\ot C)=S^a(B\ot C)\ot S^b(B\ot C)-S^{a+1}(B\ot C)\ot S^{b-1}(B\ot C)$, 
so we first compute 
$$\begin{array}{l}
\bigoplus_{a\ge b\ge 0}t^{a+b}S_{a,b}A\ot S^a(B\ot C)\ot S^b(B\ot C)= \\
\hspace*{2cm} =
\bigoplus_{\a\ge \b, \g\ge\d, \a+\b\ge\g+\d}t^{\a+\b+\g+\d}S_{\a+\b,\g+\d}A\ot 
S_{\a,\b}B\ot S_{\g,\d}
B\ot S_{\a,\b}C\ot S_{\g,\d}C.
\end{array}$$ 
The last equality follows from Cauchy formula. Now the Clebsh-Gordon formula 
implies that 
$S_{\a,\b}B\ot S_{\g,\d}B=S^{\a-\b}B\ot S^{\g-\d}B=\oplus_{0\le k\le\a-\b,\g-\d}
S^{\a-\b+\g-\d-2k}B$. 
Define the formal series $P_{u,v,w}(t)$ by the identity 
$$\bigoplus_{a\ge b\ge 0}t^{a+b}S_{a,b}A\ot S^a(B\ot C)\ot S^b(B\ot C)=
\bigoplus_{u,v,w\ge 0}P_{u,v,w}(t)S^uA\ot S^vB\ot S^wC,$$
and observe that the coefficient of $t^n$ inside  $P_{u,v,w}(t)$ is equal to the number 
of solutions 
of the system of equations in nonnegative integers
$$\left\{\begin{array}{rcl} 
n & = & \a+\b+\g+\d, \\ u & = & \a+\b-\g-\d, \\
v & = & \a+\g-\b-\d-2k, \\ w & = & \a+\g-\b-\d-2l, 
\end{array}\right.$$
with $\a+\b\ge\g+\d$ and $0\le k,l\le\a-\b,\g-\d$. From these equations we first deduce that 
$u+v=2\a-2\d-2k$ and $u+w=2\a-2\d-2l$, which imply that $u,v,w$ have the same parity. 
Let $2r=u+v$ and $2s=u+w$, so that $k=\a-\d-r$ and $l=\a-\d-s$. Suppose that $u\ge v\ge w$, 
so that in particular $r\ge s$. Then 
$$
P_{u,v,w}(t)  =  \sum_{\substack{\g+s\ge\a\ge\d+r\\ \d+s\ge\b\ge 0\\ \a+\b=\g+\d+u}}
t^{\a+\b+\g+\d}
=\sum_{\substack{\a\ge\d+r\\ \d+s\ge\b\ge 0\\ \b+s=\d+u}}t^{2\a+2\b-u} =
\frac{t^v}{1-t^2}\sum_{\substack{\d+s\ge\b\ge 0\\ \b+s=\d+u}}t^{2\d+2\b}
 = \frac{t^{u+v-w}(1-t^{2w+2})}{(1-t^2)^2(1-t^4)}.$$
A   similar computation shows that 
$$\bigoplus_{a\ge b>0}t^{a+b}S_{a,b}A\ot S^{a+1}(B\ot C)\ot S^{b-1}(B\ot C)=
\bigoplus_{u,v,w\ge 0}Q_{u,v,w}(t)S^uA\ot S^vB\ot S^wC,$$
where  $Q_{u,v,w}(t)=\frac{t^{u+v-w+2}(1-t^{2w+2})}{(1-t^2)^2(1-t^4)}$ for $u\ge v\ge w$.
Thus 
$$\bigoplus_{k\ge 0}t^k S^k(A\ot B\ot C) =\bigoplus_{u,v,w\ge 0}
\frac{t^{u+v+w-2m}(1-t^{2m+2})}{(1-t^2)(1-t^4)}S^uA\ot S^vB\ot S^wC,$$
with the notation $m=\min(u,v,w)$. The lemma is now just a transcription
 of this formula.\qed

\medskip The lemma can   be rewritten in the following form:

$$\begin{array}{l}
\bigoplus_{k\ge 0}t^k S^k(A\ot B\ot C) = \frac{1}{(1-tA\ot B\ot C)(1-t^3A\ot B\ot C)(1-t^4)}
\times \\ \hspace*{1cm}  \times\Biggl(\frac{1}{1-t^4S^2A\ot S^2B} 
\Bigl(\frac{1}{1-t^2S^2A} +\frac{1}{1-t^2S^2B}-1\Bigr) +
\frac{1}{1-t^4S^2B\ot S^2C}\Bigl(\frac{1}{1-t^2S^2B}+\frac{1}{1-t^2S^2C}-1\Bigr) +
\\ \hspace*{2cm} +
\frac{1}{1-t^4S^2C\ot S^2A}\Bigl(\frac{1}{1-t^2S^2C}+\frac{1}{1-t^2S^2A}-1\Bigr) 
-\frac{1}{1-t^2S^2A}-\frac{1}{1-t^2S^2B}-\frac{1}{1-t^2S^2C}+1\Biggr).
\end{array}$$
and the theorem follows.

\subsection{Isotropy representations of orthogonal adjoint varieties}  

The set of semi-simple parts
of the isotropy groups  for all fundamental adjoint
varieties consists of  the subexceptional
series plus   $\fsl_2 \times \fso_n$ acting
on $V=A\ot B= \BC^2\ot \BC^n$. This new case is quite similar 
as:

\begin{prop}
 
$$\bigoplus_{k\ge 0}t^k S^kV =\frac{1}{ 
(1-tV) (1-t^2S_{[1,1]}B) (1-t^3V) (1-t^4)}\Bigl(\frac{1}{1-t^2S_{[2]}A}+
\frac{1}{1-t^4S_{[2]}B}-1\Bigr).$$

\end{prop}

Thus the
covariant algebra $\CC[V]^{\fu}$ is not a polynomial 
algebra as in the subexceptional cases, 
although it has generators of exactly the same degrees.
The fact that we no longer obtain a polynomial 
algebra seems to be related, first to the nonsimplicity of $\fg=S^2A\op S_{[1,1]}B$, 
and also to the fact that $V_2=S_{[2]}B\op S^2A\ot S_{[1,1]}B$ partly comes 
from $\fg$, since its 
second factor is just the tensor product of the two components of $\fg$. 
 Unlike the subexceptional case, the orbit closures here are not nested.

\begin{proof}
The Cauchy formula   gives
$$S^k(A\otimes B)=
\bigoplus_{\substack{l\ge m\ge 0\\ l+m=k}}S_{l,m}A\otimes S_{l,m}B.
$$
Here the Schur power $S_{l,m}B$ is not irreducible as a $\fso_n$-module, its decomposition 
into irreducibles can be found in \cite{litt} and is given by
$$S_{l,m}B = \bigoplus_{\substack{a\ge b\ge 0,\\ p\ge q\ge 0}}c_{(2a,2b),(p,q)}^{l,m}S_{[p,q]}B,$$
where $S_{[p,q]}B$ denotes the irreducible $\fso_n$-module indexed by the two-parts 
partition $(p,q)$, and the Littlewood-Richardson coefficient $c_{(2a,2b),(p,q)}^{l,m}$
is the multiplicity of the $GL(C)$-module $S_{l,m}C$ inside the tensor product 
$S_{2a,2b}C\ot S_{p,q}C$, where $C$ is some vector space of dimension at least two. 
By the Littlewood-Richardson rule, this multiplicity equals the number of triples of 
non-negative integers $\a,\b,\g$ such that $0\le\b\le 2a-2b$ and $0\le\g\le\a$, $l=2a+\a$,
$m=2b+\b+\g$, $p=\a+\b$ and $q=\g$. Letting $a=b+c$,  we get
$$\bigoplus_{k\ge 0}t^k S^kV = \bigoplus_{\substack{b,c,\a,\b,\g\ge 0\\0\le\b\le 2c, 0\le\g\le\a}}
t^{4b+2c+\a+\b+\g}S_{2c+\a-\b-\g}A\ot S_{[\a+\b,\g]}B.$$
We   let $\a=\g+\d$, and for $a$ we distinguish two cases: either $\b=2\rho$ is even and we let 
$a=\rho+\s$, or  $\b=2\rho+1$ is odd and we let $a=\rho+\s+1$. Then 
$$\begin{array}{rcl}
\bigoplus_{k\ge 0}t^k S^kV & = & \frac{1}{1-t^4}\Bigr(
\bigoplus_{\rho,\s,\g,\d\ge 0}t^{4\rho+2\s+2\g+\d}S_{2\s+\d}A\ot S_{[\g+\d+2\rho,\g]}B + \\
 & & \hspace*{2cm}+\bigoplus_{\rho,\s,\g,\d\ge 0}t^{4\rho+2\s+2\g+\d+3}S_{2\s+\d+1}A\ot 
S_{[\g+\d+2\rho+1,\g]}B\Bigr),
\end{array}$$
giving the rational expressions
$$\begin{array}{rcl}
\bigoplus_{k\ge 0}t^k S^kV & = & 
\frac{1+t^3A\ot B}{(1-t^4)(1-tA\ot B)(1-t^2S_2A)(1-t^2S_{[1,1]}B)(1-t^4S_{[2]}B)} \\
 & = & \frac{1-t^6S_2A\ot S_{[2]}B}{(1-tA\ot B)(1-t^4)(1-t^3A\ot B)(1-t^2S_2A)
(1-t^2S_{[1,1]}B)(1-t^4S_{[2]}B)} \\
 & = & \frac{1}{(1-tA\ot B)(1-t^4)(1-t^3A\ot B)(1-t^2S_{[1,1]}B)}
\Bigl(\frac{1}{1-t^2S_2A}+\frac{1}{1-t^4S_{[2]}B}-1\Bigr).
\end{array}$$
\end{proof}

\bigskip

\section{The Severi series}

Zak proved Hartshorne's conjecture
that a smooth subvariety $X^n\subset\pp{n+a}$ not contained in a hyperplane
cannot have a degenerate secant variety
if $a<\frac n2+2$, and then
classified the boundary case. The answer   gives rise to the
series corresponding to the second line in Freudenthal's square:
$$A_2, \quad A_2\times A_2,\quad A_5, \quad E_6$$
which we parametrize by $m=1,2,4,8$.
We could add the finite group ${\mathfrak S}_3$ with $m=0$. (In the case $m=0$ that
$V$, defined below,  has the correct dimension, but $\fg$ does not.)

 Freudenthal's incidence geometries
 \cite{freud} or the triality model \cite{LM3}
 distinguishes
 two isomorphic representations of dimension $3m+3$. We choose one, call it $V$
 and call its dual $V^*$. In fact $V^*=V_Q=I_2(X)$ with respect to
 our previous notations, where $X\subset \BP V$ denotes
  the unique closed orbit. 
While not singled out by the triality
 model, $\fg$ does occur as $\fg = (VV^*)_{Aad}$, i.e., as a
 space of linear syzygies. Its dimension is
 $\dim\fg = \frac{4(m+1)(3m+2)}{m+4}$.

\medskip
Let $\Gamma_0$ be the
automorphism group of the Dynkin diagram, and $\Gamma\subset\Gamma_0$
be the subgroup preserving the marked Dynkin diagram ($\Gamma$ is
trivial except for $\fg=A_2\times A_2$, for which
$\Gamma={\mathfrak S}_2$). We obtain  the following
decomposition formulae into  $\fg\times\Gamma_0$-Casimir eigenspaces:

$$\begin{array}{rcl}
\Lambda^k V & = & V_k, \quad  2\le k\leq 6, \\
 \fg\ot V & = & V\op {V_2}^*\op {V\fg} \op J,  \\
S_{21} V & = & \fg\op {VV^*}\op {VV_2},  \\
S_{31}V & =& V\op {V_2}^*\op {V\up 2V^*}\op {V\fg}\op V^*V_2\op
{V\up2V_2},
\\
S_{22} V & =& V\op {V\up 2}^*\op {V\up 2V^*}\op {V\fg}\op
{V_2\up 2}, \\
 S_{211} V & \supset& {V_2}^*\op {V\fg}\op J\op V^*V_2.
\end{array}
$$

The Severi series is sub-minuscule, so Theorem 5.1  applies. There
are only three orbits:
$$\bigoplus_{k\ge 0}t^k S^kV = 
(1-tV)^{-1}(1-t^2V^*)^{-1}(1-t^3)^{-1}.$$

The   Casimir eigenvalues for these modules are 
all of the form $\frac{am+b}{9m}$ with $a,b\in \BZ$ 
and are linear functions of $(\la, \la )$ as before.

The dimensions of these modules are rational functions of $m$
with simple denominators. The formulas
can be found in  \cite{LM3}, with the exceptions of $V_k$
which is obvious, $\fg$ given above and
$$\dim J =\frac{3(m+1)(8-m)(3m+2)}{2(m+4)(m+6)}.$$

As explained in \S 2, $\fg\subset S_1$ is a subspace of the space
of linear syzygies of $X\subset\ppp V$. In fact more is true:
 
\begin{prop} 
Let   $S_k$ denote the chain of linear syzygies in the minimal
resolution of a Severi variety, beginning with $S_0=I_2(X)$. Then
$$S_0=V^*, \qquad S_1=\fg, \qquad S_2=J.$$
\end{prop}

The highest weights of the modules involved in the decomposition
formulas
are given in the following table. Note that for $A_1^{\op 2}$ each time a representation
occurs, its mirror occurs as well which we supress in the list. In particular, the adjoint
representation is not irreducible:
\medskip

$$\begin{array}{lllll}
 & {A}_2 & A_2^{\op 2} & A_5 & E_6 \\
  & & & & \\
V & \left[2,0\right] & \left[0,1|1,0\right] & \left[0,1,0,0,0\right]
& \left[1,0,0,0,0,0\right] \\
\fg & \left[1,1\right] & \left[1,1|0,0\right] &
\left[1,0,0,0,1\right] & \left[0,1,0,0,0,0\right] \\
 {V_2} & \left[2,1\right] & \left[1,0|2,0\right] &
\left[1,0,1,0,0\right] & \left[0,0,1,0,0,0\right] \\
 {V_3} & \left[3,0\right]    & \left[ 0,0 | 3,0\right] &
\left[0,0,2,0,0\right] & \left[0,0,0,1,0,0\right] \\
  & \left[0,3\right]     & &
\left[2,0,0,1,0\right] &   \\
  &     & \left[ 1,1 | 1,1\right] &
 &   \\
{V_4} & \left[1,2\right]    & \left[2,0|0,2\right] &
\left[3,0,0,0,1\right] & \left[0,1,0,0,1,0\right] \\
  &   & \left[ 1,0 | 1,2\right] &
\left[1,0,1,1,0\right] &   \\
 {V_5} & \left[0,2\right]   & \left[2,0|0,2\right]  &
\left[4,0,0,0,0\right]  & \left[0,2,0,0,0,1\right] \\
  &   & \left[ 1,0 | 1,2\right]  &
 \left[2,0,1,0,1\right] & \left[0,0,0,0,2,0\right] \\
 &   &   \left[ 2,1 | 0,1\right]  &
 \left[0,1,0,2,0\right] &  \\
{V_6} & \left[0,0\right]   &  \left[ 0,0 | 3,0\right] &
\left[1,1,0,1,1\right]  & \left[0,3,0,0,0,0\right] \\
  &   &  \left[ 0,3 | 0,0\right] &
 \left[3,0,1,0,0\right] & \left[0,1,0,0,1,1\right] \\
 &   &  \left[ 1,1 | 1,1\right]  &
 \left[0,0,0,3,0\right] &  \\
  J & \left[0,1\right] & \left[1,0|0,1\right] & \left[2,0,0,0,0\right] & \\
\end{array}$$
\bigskip

\section{The Severi-section series}

This is the series of the first line in Freudenthal's square:
$$A_1, \quad A_2,\quad C_3, \quad F_4.$$

\medskip
Again let $m=1,2,4,8$ respectively. This series does
not correspond to a line in Vogel's plane, but 
$B_1=A_1,A_2,C_3$ are on a line, their
parameters being $(7m-8,-2m,4)$ for $m=1,2,4$. In particular the
sum of these coefficients is $5m-4$, which is precisely the
denominator in the Casimir eigenvalues below.
There is a distinguished $\fg$-module $V$ of dimension $3m+2$.

\medskip

We have a uniform 
decomposition  
 $$
\La 2 V = \fg\op {V_2} 
$$
where the presence of both factors is easily understood, the
first because $\fg$ preserves a quadratic form on $V$
and thus lies in $\fso (V)$, the second by diagram induction. 

\medskip

Let $\Gamma_0$ be the automorphism group of the Dynkin diagram, and $\Gamma\subset\Gamma_0$
be the subgroup preserving the marked Dynkin diagram ($\Gamma$ is trivial except for $\fg=A_2$, 
for which $\Gamma={\mathfrak S}_2$). We obtain  the following decomposition formulae into 
irreducible $\fg\times\Gamma $-modules:

\begin{prop} Let $\e_m=1$ for $m=1$, $\e_m=0$ for $m=2,4,8$. Then
$$\sum_{k\ge 0}t^kS^kV = \frac{1-\e_mt^6V^{(3)}}{(1-tV)(1-t^2)(1-t^2V)(1-t^3)(1-t^3V_2)}.$$
\end{prop}

All the generators except $V_2$ and the quadratic form are generators of
ideals of orbits. The presence of $V_2$
can be understood as follows:  
 the polarization of the cubic invariant gives a map 
$r: S^2V\ra V^*\simeq V$, hence a map $s: S^3V\ra V_2$ by 
letting $s(v^3)=p(v,r(v))$, where we identify $V\simeq V^*$ using
the quadratic form.

\medskip
For $m=4$ or $8$ the representation $V$ is again exceptional in the sense of \cite{brion}, 
whose results imply the Proposition in those cases. 

In the case $m=1$ the invariant algebra $\CC[V]^{\fg}$ is free, but there exists a 
(unique) relation in degree six between the fundamental covariants in $\CC[V]^{\fu}$. 
This is the classical case of quartic binary forms
(see \cite{dixm} and references therein for covariants of binary forms).

\medskip
The   Casimir eigenvalues for   these modules are all
of the form $\frac{am+b}{5m-4}$ with $a,b\in \BZ$  and are
linear functions of $(\la,\la)$.

We have 
$$\dim\fg = \frac{3m(3m+2)}{m+4},\ \ \dim V_2=\frac{(3m+2) (3m+4) (m+1)}
{2(m+4)}
$$
and  highest weights  
are given in the following table:

$$\begin{array}{lrrrl}
 & A_1 & A_2 & C_3 & F_4 \\
 & & & & \\
V & \left[4\right] & \left[1,1\right] & \left[0,1,0\right] &
\left[0,0,0,1\right] \\
\fg & \left[2\right] & \left[1,1\right] & \left[2,0,0\right] &
\left[1,0,0,0\right] \\
{V_2} & \left[6\right] & \left[3,0\right] & \left[1,0,1\right] &
\left[0,0,1,0\right] \\
\end{array}$$
\medskip

\subsection{The adjoint representation of $\fsl_3$}.
The case $m=2$ deserves some explanation since $V_2$ is irreducible   as 
a $\fg\times\Gamma $-module, but has two components as a $\fg$-module: $V_2=V_{3\o_1}
\oplus V_{3\o_2}$,   the nontrivial element of $\Gamma =\ZZ/2\ZZ$ permutes the two 
components.  The identity above should be understood as 
$$\sum_{k\ge 0}t^kS^kV = \frac{1}{(1-tV)(1-t^2)(1-t^2V)(1-t^3)}
\Bigl(\frac{1}{1-t^3V_{3\o_1}}+\frac{1}{1-t^3V_{3\o_2}}-1\Bigr).$$
Note that $(1-t^3V_{3\o_1})^{-1}+(1-t^3V_{3\o_2})^{-1}-1=1+
\sum_{k>0}t^k(V_{3k\o_1}\oplus V_{3k\o_2})$, so that the preceeding identity means that 
$\fsl_3$ is exceptional in the sense that $\CC[\fsl_3]^{\fu\times\Gamma }$ is a polynomial 
algebra, although $\CC[\fsl_3]^{\fu}$ is not. 

We briefly explain how one obtains the preceeding generating function $g_{\fsl_3}(t)$ 
for the symmetric powers of $\fsl_3$. If $U$ denotes the natural three-dimensional module, 
first note that $U^*\ot U = \fsl_3\oplus\CC$, so that $g_{\fsl_3}(t)=(1-t)g_{U^*\ot U}(t)$.
Again the symmetric powers of a tensor product are given by the Cauchy formula:
$$S^k(U^*\ot U)=\sum_{a+2b+3c=k}S_{a+b+c,b+c,c}U\otimes S_{a+b+c,b+c,c}U^*.$$
But as $\fsl_3$-modules, $S_{a+b+c,b+c,c}U^* = S_{a+b,b}U^* = S_{a+b,a}U$, and we get 
$$g_{\fsl_3}(t)=\frac{1-t}{1-t^3}\sum_{a,b\ge 0}t^{a+2b}S_{a+b,a}U\ot S_{a+b,b}U.$$ 
Now we use the Littlewood-Richardson rule to compute these scalar products: we refer
the reader to \cite{man} for the statement and the terminology we use in the sequel. 
Following this rule, the irreducible components of $S_{a+b,a}U\ot S_{a+b,b}U$ are 
encoded by skew-tableaux of the following type:

\centerline{\young(\hh\hh\hh\hh\hh\hh 1111,\hh\hh 111222,1222)}

\medskip\noindent
We have $a+b$ empty boxes on the first line, $b$ on the second line. 
We add $\a_i$ boxes numbered $1$ on the $i$-th line, $i=1,2,3$, with total number $a+b$, 
and $\b_j$ boxes numbered  $2$ on the $j$-th line, $j=2,3$, with total number $a$.
Moreover, there are two types of constraints. First we must get a semistandard skew-tableau,
which means that below a box numbered $1$ there can be no box also numbered $1$, 
and below a box numbered $2$ there can be no box at all. This means that 
$$\a_2\le a, \quad\a_2+\b_2\le a+\a_1, \quad \a_3\le b, \quad \a_3+\b_3\le b+\a_2.$$
Second, the word one obtains by reading the numbered boxes right to left and top to bottom
must be Yamanouchi (or a lattice word), which means that 
$$\b_2\le\a_1 \quad {\rm and} \quad \b_2+\b_3\le \a_1+\a_2.$$
When   these conditions are fulfilled, we have   
$S_{a+b+\a_1,b+\a_2+\b_2, \a_3+\b_3}U\subset S^k(U^*\ot U)$.

Recalling that $a=\a_1+\a_2+\a_3$ and $b=\b_2+\b_3$, it is easy to see that this set
of inequalities actually reduces to:
$$\b_2\le\a_1,\quad \a_2\le\b_2+\b_3, \quad \b_2+2\b_3\le\a_1+2\a_2.$$
The first of theses implies that we can write $\a_1=\b_2+u$ for some non-negative integer $u$. 
Then we have two cases. 

If $\a_2\ge\b_3$, we let $\a_2=\b_3+v$ for some non-negative integer $v$, 
then the third inequality is automatically true and the second one reduces to $\b_2\ge v$, so that 
$\b_2=v+w$ for some non-negative integer $w$. The $\fsl_3$-module we obtain this way is 
$S_{2u+3v+2w+\a_3+\b_3, u+3v+w+\a_3+\b_3, \a_3+\b_3}U=S_{2u+3v+2w, u+3v+w}U$, and the 
overall contribution of this case is 
$$\begin{array}{l}
\sum_{u,v,w,\a_3,\b_3\ge 0}t^{2u+3v+w+2\a_3+\b_3}S_{2u+3v+2w, u+3v+w}U\\
\hspace*{8cm} =\frac{1}{(1-t^2S_{21}U)(1-t^3S_{33}U)(1-tS_{21}U)(1-t^2)(1-t^3)}.
\end{array}$$

If $\a_2\le\b_3$, we let $\a_2=\b_3-v$ for some non-negative integer $v$; 
then the second inequality is automatically true and the third one reduces to $u\ge 2v$, so that 
$u=2v+w$ for some non-negative integer $w$. The $\fsl_3$-module we obtain this way is 
$S_{4v+2w+\a_2+\a_3+2\b_2, v+w+\a_2+\a_3+\b_2, v+\a_2+\a_3}U=S_{3v+2w+2\b_2, w+\b_2}U$,
and the overall contribution of this case is 
$$\begin{array}{l}
\sum_{v,w,\a_2,\a_3,\b_2\ge 0}t^{3v+2w+2\a_2+\a_3+\b_2}S_{3v+2w+2\b_2, w+\b_2}U \\
\hspace*{8cm} =\frac{1}{(1-t^3S_{3}U)(1-t^2S_{21}U)(1-t^2)(1-t)(1-tS_{21}U)}.
\end{array}$$

Finally, we counted  the case $\a_2=\b_3$ twice, whose contribution is easily calculated to be 
$$\sum_{u,\a_2,\a_3,\b_2\ge 0}t^{2u+\a_2+2\a_3+\b_2}S_{2u+2\b_2, u+\b_2}U
=\frac{1}{(1-t^2S_{21}U)(1-t)(1-t^2)(1-tS_{21}U)}.$$
Putting together theses three contributions we easily obtain the expression we claimed for
the generating series $g_{\fsl_3}(t)$.

\bigskip

\section{The highest possible Casimir eigenspace of $\La k V$}

Let $V$ be a fundamental representation of a  simple Lie algebra $\fg$
with highest weight $\l$ and 
  Casimir eigenvalue $\th_V$. Let $\a$ denote the   simple root
whose coroot is Killing-dual to $\l$.
 Define  $V_k\subseteq \La k V$
to be the (possibly empty) subspace with Casimir eigenvalue 
$$\th_{V_k}:= k\th_V+ k(k-1)\left[ (\lambda, \lambda) - (\a,\a) \right].
$$

We expect that $V_k$, when nonempty,
 is the highest Casimir eigenspace in $\La kV$.
We show below that this is the case when $V$ is minuscule,
  it is true when $V$ is adjoint by \cite{kos1}, and we
extend it to other fundamental representations in low
degrees in proposition 9.3. Let $k_0$ denote
the largest $k$ for which $V_k$ is nonempty.

\rem In \cite{kos1}, a beautiful characterization
of $V_k$ is given in the case
$V=\fg$ is the adoint representation: the components of $\fg_k$ correspond
to abelian ideals of
a fixed Borel $\mathfrak b$. Our answer in the general case
is not as elegant and it would be nice to have a simpler characterization.

For the adjoint representations  $k_0$ is explicitly known. Also
note that $\th_{\fg_k}=k$ as our formula predicts. 
\medskip

In the standard representations of classical series and the Severi series, we
 have $V_k=\La kV$ in low degrees. In the subexceptional series,  in low degrees
$V_k$ is the primitive subspace for the symplectic form $\omega$, i.e.,
$\La k V= V_k \op (\omega\ww\La {k-2}V)$. For the exceptional series,
at least in low degrees,
$\fg_k$  is the primitive part of $\La k\fg$
  For example, $\La 2\fg = \fg \op \fg_2$,
but the inclusion $\fg\ra\La 2\fg$ is just the Lie
bracket,
so the only primitive piece is $\fg_2$.

\medskip
 Let $W_k$ denote the Casimir eigenspace of $\La k V$ of maximal 
eigenvalue. The discussion of \cite{kos1} implies that $W_k$ is decomposably 
generated, i.e.,  its highest weight vectors are all of the 
form $v_1\we\cdots \we v_k$ for some weight vectors $v_1,\ldots ,v_k$ of $V$. 
(Kostant only considered the case where $V=\fg$ is the adjoint  
representation, but his arguments   apply to any irreducible module.) 
Note that   the set of vectors $v_1\hd v_k$ is $B$-stable
and conversely a $B$-stable set of vectors wedged together
furnishes a highest weight vector. Here $B$ denotes
the Borel compatible with our choices. 

We will call a $B$-stable set of weight vectors   {\it complete}.
We will also call the corresponding set of weights complete.
A subset $S$ of the weights of $V$ is complete if and only if
for all $\mu \in S$ and each  $\b$ that  is a sum of positive roots, 
if $\mu +\b $
is a weight of $V$,
then $\mu +\b\in S$.

Thus the problem of characterizing $W_k$ is to characterize
which complete subsets of weights (possibly with multiplicities, bounded by their 
multiplicities in $V$)  determine a maximal Casimir eigenvalue.


Let $H_j$ be an orthonormal  basis of the Cartan subalgebra
of $\fg$, and $X_{\b}$ a generator of the 
root space $\fg_{\b}$. Let $\TH$ denote the Casimir operator. We have
(see \cite{LM2}) 
$$
\begin{aligned}
\TH(v_1\ww\cdots\ww v_k)&=\sum_iH_iH_i(v_1\ww\cdots \ww v_k)
+\sum_{\b\in\Delta}X_{\b}X_{-\b}(v_1\ww\cdots \ww v_k)\\
&=
k\th_V v_1\we\cdots \we v_k
+2\sum_l\sum_{i<j}v_1\we \cdots \we H_lv_i\we\cdots \we H_lv_j\we\cdots \we v_k \\
 & \hspace*{1cm}  +\sum_{\b\in\Delta}\sum_ {i\ne j}\frac{4}{(X_{\b},X_{-\b})}
v_1\we \cdots \we X_{\b}v_i\we\cdots \we X_{-\b}v_j\we\cdots \we v_k.
\end{aligned}
$$

\medskip
In order to state the main result of this section, we
define the {\it diameter} of a subset $S$ of the weights of $V$
to be the minimal number $\delta$ such that $\norm{\mu-\mu '}^2
\le \delta (\a,\a)$ for all $\mu,\mu '\in S$.
 The diameter of $V$ is obtained for 
$\mu=\l$ and $\mu '=w_0(\l)$, where $w_0$ denotes the longest element of the Weyl group. 
Thus, we easily compute that $\delta=i$ for the $i$-th fundamental representation of $A_l$, 
$\delta=2$ for the natural representations of $C_l$, $\delta =l$ for the spin 
representation of
$B_l$, $\delta = [l/2]$ for a spin representation of $D_l$, $\delta=2$ for the minuscule 
representation of $E_6$, and $\delta=3$ for that of $E_7$. Note that when $\delta=2$, any 
decomposably generated component of $\La k V$ has maximal Casimir eigenvalue.

\begin{prop} Let $V$ be a minuscule representation. Then the irreducible components 
of $\La k V$ have Casimir eigenvalue less than or equal to $\th_{V_k}$. Those 
with Casimir eigenvalue equal to $\th_{V_k}$ are in correspondance with complete
cardinality $k$ subsets $S$ of the set of weights of $V$ 
of diameter at most $2$. In the case of the minuscule representation of $B_l$,
we require additionally   that the difference between
two   elements  of $S$ cannot be a root strictly longer than $\a_l$. 
\end{prop}

\proof 
Let $U$ denote a component of $\La k V$ of maximal Casimir eigenvalue, let 
$v_1\ww\cdots 
\ww v_k$ be a highest weight vector. 
Let  $\mu_i$ denote the weight of $v_i$, and suppose that $V$ is endowed 
with an invariant  Hermitian product $\langle\ ,\ \rangle$, such that the $v_i$ are part 
of a unitary basis.  Then the eigenvalue of the Casimir operator on $U$ is 
$$
\th_U
=\langle\Theta (v_1\we\cdots \we v_k),v_1\we\cdots \we v_k\rangle =
 k\th_V+\sum_{i\ne j}(\mu_i,\mu_j) +\sum_{\b\in\Delta}\sum_ {i\ne j}
\frac{\langle X_{\b}v_i\we X_{-\b}v_j,v_i\we v_j\rangle}{(X_{\b},X_{-\b})}.
$$
Since weight vectors of distinct weights are orthogonal, $\langle X_{\b}v_i
\we X_{-\b}v_j,v_i\we v_j\rangle$ can be nonzero only if $\mu_i=\mu_j-\b$ and there 
exist  scalars $s$ and $t$ such that $X_{\b}v_i=sv_j$ and $X_{-\b}v_j=tv_i$. 
Assuming this, we compute 
$$stv_i=sX_{-\b}v_j=X_{-\b}X_{\b}v_i=[X_{-\b},X_{\b}]v_i+X_{\b}X_{-\b}v_i.$$
The latter term is zero since, $V$ being minuscule, $X_{-\b}^2v_j=0$ (\cite{bou}, page 128).
Moreover, we may suppose that $[X_{-\b},X_{\b}]=H_{\b}$ is
the coroot of $\b$ (see
\cite{bou}, page 82), and note that in this case $2(X_{\b},X_{-\b})=-(H_{\b},H_{\b})
$, and we get 
$$\frac{\langle X_{\b}v_i\we\cdots \we X_{-\b}v_j,v_i\we v_j\rangle}{(X_{\b},X_{-\b})}
=\frac{2}{(H_{\b},H_{\b})}\mu_i(H_{\b})=(\mu_i,\b)=(\mu_i,\mu_j-\mu_i).$$ Hence 
$$\begin{array}{rcl}\th_U  & = & 
 k\th_V+\sum_{i\ne j}(\mu_i,\mu_j) +\sum_{\mu_j-\mu_i\in\Delta}(\mu_i,\mu_j-\mu_i) \\
 & = & k\th_V+\sum_{\mu_j-\mu_i\notin\Delta}(\mu_i,\mu_j) 
+\sum_{\mu_j-\mu_i\in\Delta}(\mu_i,2\mu_j-\mu_i).
\end{array}$$
Note that since $V$ is minuscule, the weights $\mu_i$ are all conjugate under the Weyl 
group, in particular they have the same norm as $\l$. We need the following   
observation: 

\begin{lemm} For $i\neq j$, either
$\norm{\mu_i-\mu_j}^2=(\a,\a)$ and $\mu_i-\mu_j\in\Delta$,
or $\norm{\mu_i-\mu_j}^2\ge 2(\a,\a)$ and
$\mu_i-\mu_j\notin\Delta$.
\end{lemm}

\proof We may suppose that $\mu_i=\l$, since the Weyl group acts transitively on the 
weights of $V$. Since $\l$ is the highest weight of $V$, we can write $\mu_j=\l-\sum_kn_k\a_k$
for some non-negative integers $n_k$, where the $\a_k$ are the simple roots. Since $\l$ is  
fundamental it is orthogonal to every simple root except $\a=\a_l$, say, and we get 
$(\mu_i,\mu_j)=(\l,\l)-n_l(\a_l,\o_l)=(\l,\l)-n_l(\a,\a)/2$, hence $\norm{\mu_i-\mu_j}^2=
n_l(\a,\a)$. 

Suppose that $n_l=1$. The highest weight of $V$ after $\l$ is $\l-\a$.
Since every nonzero weight of $V$ is obtained by a sequence
of simple reflections in $\l$, 
  there is a 
sequence $\n_i$,
$1\leq i\leq k+1$ of weights of $V$ such that $\n_0=\l$, $\n_1=s_{\a}(\l)=\l-\a$, 
$\n_t=\mu_j$ for some $t$ and $\n_{k+1}=s_{\b_k}(\n_k)$ for some simple root $\b_k$, 
which is different from $\a$ if $k\ne 0$
because $n_l=1$. But then $s_{\b_k}(\l-\n_k)=\l-\n_{k+1}$, 
thus $\l-\n_k$ is a root if and only if  $\l-\n_{k+1}$ is 
also a root. Since $\l-\n_1=\a$ is 
indeed a root, we conclude that $\mu_i-\mu_j = \l-\n_t$ is a root. 
This argument is reversible, proving the lemma if we remember the formula 
$\norm{\mu_i-\mu_j}^2=n_l(\a,\a)$. \qed 

\medskip
To conclude the proof of the proposition, we just need, for each pair $\mu_i,\mu_j$,
to choose an element $w$ of the Weyl group such that $w(\mu_i)=\l$, and define the 
integer $n_{i,j}$ to be the coefficient of $\l-w(\mu_j)$ on the simple root $\a$. 
Then $(\mu_i,\mu_j)=(\l,\l)-n_{i,j}(\a,\a)/2$ and we get the formula
$$\th_U = k\th_V+\sum_{\mu_j-\mu_i\notin\Delta}((\l,\l)-n_{i,j}(\a,\a)/2) 
+\sum_{\mu_j-\mu_i\in\Delta}((\l,\l)-n_{i,j}(\a,\a)).$$
The $n_{i,j}$ are all positive, and they are at least equal to two in the first 
sum. We conclude that $\langle\Theta (v_1\we\cdots \we v_k),v_1\we\cdots \we v_k\rangle$
will be maximal when $n_{i,j}$ is always equal to two in the first sum,
meaning that 
two weights whose difference is not a root have the square of their distance equal 
to $2(\a,\a)$, and always equal to one in the second sum (which means that their difference
is a root which is not longer than $\a$). Then we get 
$\th_U=k\theta_V+k(k-1)((\l,\l)-(\a,\a))$,
and the proposition is proved. \qed 

\medskip
In general, it is clear from the Proposition that $V_k$ is nonzero when $k$ is not too big,
but the maximal integer $k_0$ for which this is true is not so easy to compute. At least 
can we say that $k_0$ can be quite large. Indeed, for the $i$-th fundamental representation 
of $A_l$, the set of weights $\mu_{j,k}=\o_i-\e_j+\e_k$, where $1\le j\le i$ and 
$i<k\le l+1$, form, with $\o_i$, a set of weights with the required properties, 
so that $k_0>i(l+1-i)$.  We suspect that $k_0=i(l+1-i)+1$ in that case but we have not 
proved it. Note also that the number of irreducible components in $V_k$ can be arbitrary 
large, as easily follows from   Proposition 9.1. 

\medskip
A nice consequence of the fact that Proposition 9.1 above holds for the fundamental 
representations of $A_l$ is that we can   extend its validity as follows:

\begin{prop} Let $V$ be a fundamental representation of the simple Lie algebra $\fg$. 
Suppose that the corresponding node of the Dynkin diagram of $\fg$ is on an
 $A_l$-chain in $D(\fg)$, 
at distance at least $k_1$ from an extremity of the
diagram. Then for $k\le k_1+2$, $\th_{V_k}$ is 
the largest Casimir eigenvalue of $\La k V$, and the irreducible components of $V_k$ 
can be described in exactly the same way as in the preceeding proposition. \end{prop}

\proof An irreducible component of $\La k V$ with maximal Casimir eigenvalue is 
decomposably generated, hence generated by the wedge product of weight vectors whose set 
of weights form a complete subset of the set of weights of $V$.
 Moreover, there are at 
most $k$ distinct weights in this set (possibly less if $V$ has weights with multiplicity 
greater than one). But for $k\le k_1+2$, every weight of a complete $k$-set of weights 
of $V$ is of the form $\l-\theta$, where $\theta$ is a sum of simple roots corresponding to 
nodes on the $A_l$-chain only, and such that $\l-\theta$ is also a weight of the 
corresponding 
fundamental representation of $A_l$.  Indeed, we know that the weights of $V$ are the 
weights of the convex hull of the translates of $\l$ by the Weyl group, which are congruent 
to $\l$ modulo the root lattice. We can obtain the translates of $ \l $ by applying 
successively 
the simple reflections of the Weyl group so that the distance to $\l$ increases (if we measure
that distance by the sum of the coefficients of the difference, expressed in terms of simple 
roots). At the beginning of this process, the simple reflections involved are those 
associated  
to nodes of the $A_l$-chain only, and the weights one obtains are formally the same as for 
the corresponding fundamental 
representation of $A_l$. More precisely, this is the case 
until  we do not apply more than $k_1+1$ simple reflections. Moreover, we obtain no new weight 
by considering  the convex hull of those, and we conclude that the weights of $V$,
at a distance at most $k_1+1$ from $\l$, are formally the same as those of the corresponding 
representation of $A_l$, with the same multiplicity, one. The scalar products of two 
such weights can be computed in terms of 
$(\l,\l)$ and $(\a,\a)$, and a part of the Cartan matrix which only involves the $A_l$-chain,
thus the computation is formally the same as in the weight lattice of $A_l$, and therefore
the computation of the Casimir eigenvalue of a $k$-set of weights will again be formally 
identical. Finally, since we only need to consider the same $k$-sets of weights and the same 
Casimir eigenvalues  as in the $A_l$-case, the conclusions of the preceeding proposition 
for the fundamental representations of $A_l$ directly apply to $V$, and this 
concludes the proof. \qed

\medskip Returning to geometry, we arrive at the following statement, which could 
also be deduced from \cite{LM1}. 

\begin{coro} Let $V$ be a fundamental representation of $\fg$,
and let $X\subset\ppp V$ be the closed orbit. If the Fano variety $ \BF_k(X) $
of $\PP^{k-1}$'s in $X$ is nonempty, then its linear span $\langle\BF_k(X)\rangle$
is contained in $V_k$ and in this case $V_k$
is the highest Casimir eigenspace.\end{coro} 
 
\proof We know from \cite{LM0} that the closed orbits in $ \BF_k(X) $ are in correspondance 
with marked subdiagrams of type $(\fa_{k-1},\o_1)$. By the proposition above, such
subdiagrams detect components of $V_k\subset \La kV$. \qed

\medskip
We can be more precise for $k=2$:

\begin{coro} Let $V$ be a fundamental representation of $\fg$.  Then
$V_2$ is irreducible and coincides with $\langle\BF_2(X)\rangle$, the linear 
span in $\La 2V$ of the set of lines contained in the closed orbit of $\PP V$. 
\end{coro}

\end{document}